\documentclass[11pt]{amsart}
\usepackage{amssymb}
\usepackage{amsfonts}
\usepackage{amsmath}
\usepackage{amsthm}
\usepackage{dsfont}
\usepackage{mathtools}
\usepackage{geometry}
\usepackage{caption}
\usepackage{algorithm} 
\usepackage{algpseudocode} 
\usepackage{amssymb} 
\usepackage[colorlinks=true,linkcolor=blue,urlcolor=blue]{hyperref}
\mathtoolsset{showonlyrefs}
\usepackage{natbib}
\usepackage{float}

\def\d{\mathrm{d}}

\def\P{{\mathsf P}}

\def\E{{\mathsf E}}

\def\m{n}
\def\stmh{M^*}

\renewcommand{\mid}{\,|\,}

\newtheorem{theorem}{Theorem}
\newtheorem{lemma}[theorem]{Lemma}

\newtheorem{corollary}[theorem]{Corollary}
\newtheorem{definition}{Definition}
\newtheorem{assumption}{Assumption}
\theoremstyle{definition} 
\newtheorem{remark}{Remark}

\def\cX{\mathcal{X}}

\mathtoolsset{showonlyrefs} 
\usepackage{float}
\usepackage{caption}
 
\numberwithin{equation}{section}
\usepackage{enumerate}

\title[Restricted Spectral Gap Decomposition for Simulated 
Tempering]{Restricted Spectral Gap Decomposition for Simulated \\
Tempering Targeting Mixture Distributions}
\author[J.\ Garg, K.\  Balasubramanian ,Q.\ Zhou]{Jhanvi Garg \and Krishna Balasubramanian \and Quan Zhou$^\dagger$} 
\address{J.\ Garg: Department of Statistics, Texas A\&M University, College Station, TX 77843.}
\email{\href{mailto:gargjhanvi@stat.tamu.edu}{gargjhanvi@stat.tamu.edu}}
\address{K.\  Balasubramanian: Department of Statistics, University of California, Davis, CA 95616}
\email{\href{kbala@ucdavis.edu}{kbala@ucdavis.edu}}
\address{Q.\ Zhou: Department of Statistics, Texas A\&M University,  College Station, TX 77843.}
\email{\href{mailto:quan@stat.tamu.edu}{quan@stat.tamu.edu}}

\date{}
\numberwithin{equation}{section}

\usepackage{enumerate}
\begin{document}
\begin{abstract} 
Simulated tempering is a widely used strategy for sampling from multimodal distributions. In this paper, we consider simulated tempering combined with an arbitrary local Markov chain Monte Carlo sampler and present a new decomposition theorem that provides a lower bound on the 
restricted spectral gap of the algorithm for sampling from mixture distributions.  By working with the restricted spectral gap, the applicability of our results is extended to broader settings such as when the usual spectral gap is difficult to bound or becomes degenerate. We demonstrate the application of our theoretical results by analyzing simulated tempering combined with random walk Metropolis--Hastings for sampling from mixtures of Gaussian distributions. Our complexity bound scales polynomially with the separation between modes, logarithmically with \(1/\varepsilon\), where \(\varepsilon\) denotes the target accuracy in total variation distance, and exponentially with the dimension \(d\).
\end{abstract}

\maketitle

\begingroup
\renewcommand{\thefootnote}{\fnsymbol{footnote}}	 
  \footnotetext[2]{Corresponding author.  E-mail: \texttt{quan@stat.tamu.edu}}
\endgroup

\section{Introduction}\label{sec:intro}

Efficient sampling from complex distributions is a foundational problem with numerous applications across various fields, including computational statistics \citep{robert1999monte, liu2001monte,  brooks2011handbook, owen2013monte}, Bayesian inference \citep{gelman2013bayesian}, statistical physics \citep{newman1999monte, landau2021guide}, and finance \citep{dagpunar2007simulation}. These distributions are often multimodal, reflecting underlying heterogeneity in the data. Sampling from such distributions presents challenges closely related to those encountered in non-convex optimization, where objective functions with multiple local minima require methods capable of effectively exploring the solution space. While discretizations of Langevin dynamics (for a comprehensive overview, see \citet{chewi2024log})  excel in log-concave settings where gradients reliably guide the sampler toward a single mode, they tend to be less effective in multimodal landscapes, which require strategies capable of navigating between separated modes.  The assumption of dissipativity, which limits the growth rate of the potential, has been widely used in previous works to establish better convergence rates for Langevin Monte Carlo (LMC) in such settings \citep{raginsky2017non,durmus2017nonasymptotic, erdogdu2018global, erdogdu2021convergence, mou2022improved, mousavi2023towards}. More recently, \citet{balasubramanian2022towards} characterized the performance of averaged LMC for target densities that are only H\"{o}lder continuous, without relying on functional inequalities or curvature-based assumptions; they measured the convergence rate using the weaker Fisher information metric.

For distributions that deviate significantly from log-concavity and exhibit numerous deep modes, additional techniques are often required to ensure efficient sampling. A comprehensive discussion of the fundamental challenges in sampling from multimodal distributions, as well as an overview of major types of Markov chain Monte Carlo (MCMC) algorithms designed for this purpose---including parallel tempering, mode jumping, and Wang--Landau methods---can be found in \citet{latuszynski2025mcmc}. The recent work of \citet{koehler2024efficiently} addresses the problem of sampling from a multimodal distribution by initializing the sampler using a small number of stationary samples. They show that, under the assumption of a $k$-th order spectral gap,  the chain efficiently yields an $\varepsilon$-accurate sample in total variation (TV) distance with $\widetilde{O}( k / \varepsilon^{2})$ samples, and they further extend the result to scenarios where the score (i.e., the drift term of LMC) is estimated. \citet{lee2024convergence} studied the {sequential Monte Carlo (SMC) method~\citep{liu1998sequential, del2006sequential, chopin2020introduction, syed2024optimised}} and derived the complexity result for mixture target distributions. Additionally, denoising-diffusion-based samplers have been proposed for sampling from non-log-concave targets without relying on isoperimetric inequalities \citep{huang2023reverse, huang2024faster, he2024zeroth}. These methods reverse the Ornstein–Uhlenbeck process and require estimating score functions (i.e., gradients of the log-density) via importance sampling, which becomes particularly challenging in high-dimensional settings.

A widely used strategy for tackling multimodality is annealing or tempering, which leverages a sequence of distributions to gradually explore complex landscapes.  \citet{guo2024provable} provided a non-asymptotic analysis of annealed Langevin Monte Carlo~\citep{neal2001annealed}, highlighting its provable benefits for non-log-concave sampling. They demonstrated that a simple annealed Langevin Monte Carlo algorithm achieves \(\varepsilon^2\) accuracy in Kullback-Leibler divergence with an oracle complexity of 
$\widetilde{O}\left( d / \varepsilon^{6} \right)$. \citet{chehab2024provable} develops a comprehensive theoretical framework for tempered Langevin dynamics, providing convergence guarantees in Kullback-Leibler divergence across a variety of tempering schedules. 
Simulated tempering \citep{marinari1992simulated} has also been introduced as a method to promote transitions between modes.  By dynamically adjusting the ``temperature" of the distribution, simulated tempering effectively smooths the energy landscape, allowing the sampler to escape local modes and traverse between high-probability regions. The sampler gradually returns to the original distribution as the temperature decreases.  A version of simulated tempering was also used in \citet{koehler2022sampling} to sample from multimodal Ising models. Parallel tempering \citep{swendsen1986replica} shares a similar mechanism with simulated tempering, but runs multiple chains in parallel at different temperatures, allowing exchanges between them to improve mixing. \citet{zheng2003swapping} established that the spectral gap of simulated tempering chain is bounded below by a multiple of the spectral gap of parallel tempering chain and a bound depending on the overlap between distributions at adjacent temperatures. \citet{woodard2009conditions} provided a lower bound on the spectral gap for a simulated tempering chain using the state space decomposition technique~\citep{madras2002markov}. This method analyzes the probability flow by dividing the state space into subsets and studying transitions within and between them. However, determining an optimal partition into subsets often requires complex spectral partitioning arguments, and estimating the spectral gap through conductance methods can introduce a squared-factor loss due to Cheeger’s inequality.

An alternative approach is to decompose the Markov chain directly rather than the state space, particularly when the target distribution is a mixture distribution or closely resembles one. This method can yield potentially tighter spectral gap bounds by separately analyzing two types of chains: local chains that efficiently explore each mixture component at every temperature level and a projected chain that governs transitions between different components and temperature levels.   This decomposition technique was introduced in \citet{ge2018simulated} as a framework for bounding the spectral gap  of the simulated tempering chain combined with LMC for mixtures of strongly log-concave distributions that are translates of each other. They established that the runtime required to reach an \(\varepsilon\) TV distance from the target distribution depends polynomially on the dimension \(d\), the mode separation \(D\), and the inverse accuracy \(1/\varepsilon\). 
 
We introduce a new Markov chain decomposition theorem for discrete-time chains, in contrast to the continuous-time framework adopted by \citet{ge2018simulated}.  
Since spectral gap bounds of discrete-Markov chains (e.g. Metropolis--Hastings algorithms) 
over unbounded Euclidean spaces is often difficult to obtain, we propose to directly decompose the restricted spectral gap~\citep{atchade2021approximate}, which, roughly speaking, can be thought of as the spectral gap of the Markov chain restricted to a subset of the state space. Geometric tools such as the path methods developed by~\citet{yuen2000applications} can be used to lower bound the restricted spectral gap. 
Intuitively, if each local chain mixes fast in a large subset of the space, where the target distribution concentrates, and the projected chain also mixes well,   
then the overall simulated tempering chain should converge fast. The projected chain can be constructed in many ways. We present a construction which differs substantially from that of \citet{ge2018simulated}; it arises naturally from the structure of the algorithm and significantly simplifies the proof. 

The remainder of the paper is organized as follows. In Section~\ref{sec:prob}, we present our decomposition theorem. In Section~\ref{sec: example}, we apply this theorem to analyze the simulated tempering combined with the random walk Metropolis--Hastings (STMH)  algorithm for sampling from mixtures of Gaussian distributions. In Section~\ref{sec:simulation}, we empirically validate our convergence guarantee for the STMH algorithm by sampling from a two-dimensional Gaussian mixture. Finally, Section~\ref{sec:conclusion} summarizes our key contributions and highlights promising directions for future research. Detailed proofs are provided in the appendix.

\section{A New Decomposition Theorem}\label{sec:prob}

\subsection{Notation and Definitions} 
We begin by introducing the necessary notation that will be used throughout the paper. We adopt the convention that uppercase letters denote probability distributions  (or transition kernels) and the corresponding lowercase letters their densities (or transition densities). For example, $P$ will be used to denote the probability distribution with density $p$.  
We use $\mathcal{L}^2(\Omega, \Pi)$ to denote the space of all real-valued functions defined on $\Omega$ that are square-integrable with respect to a measure $\Pi$. 

\begin{definition}[Restricted Spectral Gap]
\label{def:restrictedSpecGap}
Let \( K \) be a Markov transition kernel with state space \( \Omega \) and stationary probability measure \( \Pi  \). Let \( \Omega^0 \subseteq \Omega \) be measurable such that \( \Pi(\Omega^0) > 0 \).  
The \(\Omega^0 \)-restricted spectral gap of \(K\), denoted by \(\mathrm{SpecGap}_{\Omega^0 }(K)\), is defined as
\begin{equation}
    \mathrm{SpecGap}_{\Omega^0}(K)
= \inf_{g \in \mathcal{L}^2(\Omega, \Pi)} \frac{\mathcal{E}_{\Omega^0}(g, g; \Pi, K) }{ \mathrm{Var}_{\Pi, \Omega^0}(g) }, 
\end{equation}
where 
\begin{align*}
  \mathcal{E}_{\Omega^0}(g, g; \Pi, K) &=  \frac{1}{2} \int_{\Omega^0 }\!\int_{\Omega^0} \bigl[g(\omega')-g(\omega)\bigr]^2 \, \Pi(\d \omega) \,K(\omega, \d \omega'), \\
  \mathrm{Var}_{\Pi, \Omega^0}(g) &= \frac{1}{2}\int_{\Omega^0 }\!\int_{\Omega^0 } \bigl[g(\omega')-g(\omega)\bigr]^2\, \Pi(\d \omega)\,\Pi(\d \omega'). 
\end{align*} 
\end{definition} 
We will refer to $\mathcal{E}_{\Omega^0}(g, g; \Pi, K)$ as the $\Omega^0$-restricted Dirichlet form and omit $\Pi, K$ when they are clear from the context. 
When $\Omega^0 = \Omega$, $\mathrm{SpecGap}_{\Omega^0}(K)$ is known as the spectral gap of $K$, and we simply write $\mathrm{SpecGap}(K) = \mathrm{SpecGap}_\Omega(K), 
\mathcal{E}(g, g) = \mathcal{E}_{\Omega}(g, g)$ and $\mathrm{Var}_{\Pi}(g) = \mathrm{Var}_{\Pi, \Omega_0}(g)$. Note that $\mathrm{Var}_{\Pi}(g)$ equals the variance of $g(\omega)$ with $\omega \sim \Pi$.   {Intuitively, the spectral gap quantifies how rapidly a Markov chain mixes: a larger gap corresponds to faster convergence to stationarity distribution. The restricted spectral gap generalizes this idea by measuring the rate of mixing in a subset $\Omega^0$ of the state space.}

{Simulated tempering can be viewed as a bivariate Markov chain. The key idea is to augment the state space with a temperature index that allows the chain to move between a family of distributions, ranging from the target distribution (low temperature) to flattened versions (high temperature). At higher temperatures the landscape is smoother, enabling easier transitions between modes, while samples from the target are obtained by collecting states only at the lowest temperature. Formally, it is defined as follows.}

\begin{definition}[Simulated Tempering Markov Chain] \label{def: simulated}
Let $L \geq 2$ be an integer and $[L] = \{1, 2, \dots, L\}$. Let 
\( \lambda \in (0, 1) \) and \(  (r_i)_{i \in [L]} \)  be constants such that  $r_i > 0$ for each $i$, and $\sum_{i=1}^L r_i = 1.$  
Let \( M_i \), for each \( i \in [L] \), be the transition kernel of a Markov chain  with state space \( \cX \) and stationary density \( p_i \). 
The simulated tempering Markov chain has state space \( [L] \times \cX\). 
Denote its transition kernel by 
\begin{equation*}
    M = M\left( (p_i)_{i=1}^L, (r_i)_{i=1}^L, (M_i)_{i=1}^L, \lambda \right), 
\end{equation*}
which has density 
\begin{equation*}
    m( (i, x), (i',  x') ) = \left\{\begin{array}{cc}
   (1 - \lambda) m_i(x,  x'),   &   \text{ if } i = i', x \neq x', \\[5pt]
    \frac{\lambda}{2}  a( (i, x), (i', x)),   &  \text{ if } i' = i \pm 1, x = x', \\
    \end{array} \right.
\end{equation*}
where 
\begin{equation}       
    a( (i, x), (i', x)) = \min \left\{ \frac{r_{i'} p_{i'}(x)}{r_i p_i(x)}, 1 \right\}.  \label{eq:accp}
\end{equation} 
\end{definition}

In words,  the simulated tempering Markov chain evolves as follows. 
Given current state  \( (i, x) \), we sample \( u \sim \text{Bernoulli}(\lambda) \). If \( u = 0 \), draw $x'$ from $M_i(x, \cdot)$ and move to $(i, x')$. 
If \( u = 1 \),  propose \( i' = i \pm 1 \), each with equal probability \( 1/2 \), and accept the proposal with probability  $ a( (i, x), (i', x))$  (if $i' \notin [L]$, the proposal is always rejected). 
It is easy to check that the stationary distribution $P$ of the simulated tempering Markov chain  $M$  has density given by  $$p(i, x) = r_i p_i(x).$$  
 
\subsection{Decomposition of the Simulated Tempering Markov chain}

Consider the simulated tempering Markov chain given in Definition~\ref{def: simulated}. 
Let $\cX^0 \subset \cX$ be a measurable subset of $\cX$, and our goal is to derive a decomposition theorem for the spectral gap of $M$ restricted to $[L] \times \cX^0$. We assume that the  distributions $(P_i)_{i \in [L]}$ satisfies the following assumption. 

\begin{assumption}\label{ass:pi}
For each $i \in [L]$, the stationary density $p_i$ can be expressed as a mixture of $\m$ component densities: 
   \begin{equation*}
        p_i(x) = \sum_{j=1}^{\m} w_{(i,j)} p_{(i,j)}(x),
    \end{equation*}
    where \(w_{(i,j)} \geq 0\) for each $j \in [\m]$ and $\sum_{j = 1}^{\m} w_{(i,j)} = 1$.
\end{assumption}
{For simplicity, in Assumption~\ref{ass:pi} we assume that each stationary density \( p_i \) can be expressed as a mixture of simpler component densities, an assumption also made in \citet{ge2018simulated}.} Under Assumption~\ref{ass:pi}, we can introduce a latent variable $J$ so that given $I=i, J=j$, $X$ is drawn from the component distribution $P_{(i, j)}$. 
The density of the target distribution is augmented to  
    $p(i, j, x) = r_i w_{(i, j)} p_{(i, j)}(x).$ 

The overall strategy of our decomposition is  similar to Theorem~6.3 of \citet{ge2018simulated} in that we directly decompose the Dirichlet form of $M$ into the Dirichlet forms of Markov chains $(M_{(i, j)})_{i\in[L], j\in[\m]}$, where $M_{(i, j)}$ has stationary distribution $P_{(i, j)}$. 
We assume that these chains are chosen such that the following assumption is satisfied. 

\begin{assumption}\label{ass:Mij}
Let $\mathcal{E}_{i, \cX^0}$ denote  the $\cX^0$-restricted Dirichlet form of  $M_i$.  
For each $i \in [L], j \in [\m]$,  let \(  \mathcal{E}_{(i,j), \cX^0} \) be the $\cX^0$-restricted Dirichlet form of a Markov chain \( M_{(i, j)} \) on \(\cX\) with stationary density $p_{(i, j)}$. 
For each $i$ and any $g_i \in \mathcal{L}^2( \cX, P_i)$, 
    \begin{equation}\label{eq:decomp-c1}        
    \sum_{j=1}^{\m} w_{(i, j)} \mathcal{E}_{(i,j), \cX^0}(g_i, g_i) \leq C_1   \mathcal{E}_{i, \cX^0}(g_i, g_i), 
    \end{equation} 
where $C_1 > 0$ is some constant. 
Further,  each $M_{(i, j)}$ satisfies the following inequality with some constant $C_2 > 0$: 
\begin{equation}\label{eq:poincare-Mij}
     \operatorname{Var}_{P_{(i,j)}, \cX^0}(g_i) \leq C_2 \mathcal{E}_{(i,j), \cX^0}(g_i, g_i).  
\end{equation}  
\end{assumption}

{Assumption~\ref{ass:Mij} assumes a lower bound on the mixing rate of the local chains \( M_{(i, j)} \), each targeting a single component \( p_{(i, j)} \). Efficient mixing within each component is essential, as poor local mixing can limit global exploration and slow the overall convergence of the simulated tempering chain. }

\begin{remark}\label{rmk:diff-1}
    In~\citet{ge2018simulated}, each $M_i$ is a continuous-time Langevin diffusion with stationary distribution $P_i$. In this case, condition~\eqref{eq:decomp-c1} can be easily satisfied by letting 
    $M_{(i,j)}$ be the Langevin diffusion having stationary distribution $P_{(i, j)}$, since  
    \[
        \mathcal{E}_{i}(g_i, g_i) = \int \|\nabla g_i\|^2 p_i \, \d x = \int \|\nabla g_i\|^2 \sum_{j = 1}^\m w_{(i,j)} p_{(i,j)} \, \d x  = 
        \sum_{j = 1}^\m w_{(i,j)}  \mathcal{E}_{(i,j)}(g_i, g_i), 
    \]
    which yields Equation~\eqref{eq:decomp-c1} with $C_1 = 1$. However, in our setting each $M_i$ is a discrete-time Markov chain (e.g. a Metropolis--Hastings algorithm), and finding Markov chains $M_{(i, j)}$ satisfying~\eqref{eq:decomp-c1} may not be trivial.  
\end{remark}

\begin{remark}\label{rmk:diff-2}
Condition~\eqref{eq:poincare-Mij} implies that the $\cX^0$-restricted spectral gap of $M_{(i, j)}$ is at least $C_2^{-1}$. 
In contrast, \citet{ge2018simulated} requires each ``local chain'' to have a positive spectral gap.  By weakening this condition, we can develop a decomposition theorem that is particularly useful in settings where the usual spectral gap is difficult to bound over the entire state space. 
\end{remark}

Loosely speaking, condition~\eqref{eq:decomp-c1} allows us to view the dynamics of $M_i$ as governed by the hidden variable $J$. Conditional on $J = j$, the behavior of $M_i$ can be approximated by   $M_{(i, j)}$. The constant $C_2$ in~\eqref{eq:poincare-Mij} measures the convergence rate of each $M_{(i, j)}$ on $\cX^0$. 
To bound the convergence rate of the simulated tempering Markov chain $M$, we need one more assumption characterizing the  transitions between any  $(i, j)$ and $(i', j')$. 
To this end, we construct a projected chain as follows.  

\begin{definition}\label{def:project-chain}
Let  $M$ be the simulated tempering Markov chain given in Definition~\ref{def: simulated} and  Assumption~\ref{ass:pi} hold. 
Define 
\begin{equation*} 
a( (i, j, x), (i', j', x) )   =  \min \left\{ \frac{r_{i^{\prime}} w_{(i^{\prime} ,j^{\prime} )} p_{(i^{\prime} ,j^{\prime})}(x)}{r_i w_{(i,j)} p_{(i,j)}(x) }, \; 1 \right\}, \quad p_{i}( j | x) = \frac{ w_{(i, j)} p_{(i, j)}(x) }{ p_i(x)}. 
\end{equation*}
Define the projected chain with transition matrix $\overline{M}$ on $[L] \times [\m]$ by 
\begin{align*}
\overline{M}((i,j), (i^{\prime},j^{\prime})) = 
    \begin{cases}
        (1-\lambda) \displaystyle\int_{\cX^0} \frac{ p_{(i,j)}(x) }{ P_{(i, j)} (\cX^0) } \, p_i( j' | x ) \d x, & \text{if } i = i', \; j \neq j', \\[8pt] 
        \frac{\lambda}{2} \displaystyle\int_{\cX^0} \frac{ p_{(i,j)}(x) }{ P_{(i, j)} (\cX^0) } a( (i, j, x), (i', j, x) )  \d x , & \text{if } i' = i \pm 1,\;   j = j', \\[8pt] 
        1 - \sum\nolimits_{(k, l) \neq (i,j)} \overline{M}((i,j), (k,l)), & \text{if } (i',j') = (i,j). 
    \end{cases}
\end{align*} 
\end{definition}
It is easy to prove that $\overline{M}$ is indeed a transition rate matrix (i.e., all entries are non-negative and each row sums to one.) 
By checking the detailed balance condition, one obtains the following result. 
\begin{lemma}\label{lm:barM}
The stationary distribution of $\overline{M}$ is given by  
\begin{equation}\label{eq:phi-ij}
      \overline{P}(i, j)  = \frac{ r_i w_{(i,j)} P_{(i, j)}(\cX^0) }{  P([L] \times \cX^0)}. 
\end{equation}
\end{lemma}
\begin{proof}
See Appendix~\ref{sec:proofddt}.
\end{proof}

\begin{assumption}\label{ass:barM}
    Let  \(\overline{\mathcal{E}}\) be the   Dirichlet form of \( \overline{M}\).   
Then, for any $\overline{g} \colon [L] \times [\m] \rightarrow \mathbb{R}$, 
\( \overline{M} \) satisfies the following inequality for some constant $C_3 > 0$: 
\begin{equation}\label{eq:poincare-barM} 
\operatorname{Var}_{\overline{P}}(\overline{g}) \leq C_3 \mathcal{\overline{E}}(\overline{g}, \overline{g}).
\end{equation}  
\end{assumption}
{Assumption~\ref{ass:barM} assumes a lower bound on the mixing rate of the projected chain. Intuitively, if the components are well-separated, transitions between them becomes rare, resulting in slow mixing of the projected chain and, consequently, slow overall convergence of the simulated tempering chain. Conversely, when components are closer and transitions are more likely, the chain mixes more rapidly.}
\begin{remark}\label{rmk:barM}
Our construction of the projected chain $\overline{M}$ is significantly different from that in~\citet{ge2018simulated}. In particular, \citet{ge2018simulated} defined the transition between $j, j'$ by 
\begin{equation*}
\overline{M}((i,j), (i,j')) \propto \displaystyle\frac{ w_{(i,j')}}{\chi^2_{\max}(p_{(i,j)} \| p_{(i,j')})}
\end{equation*}
where $\chi^2_{\max}(P \| Q) := \max\{\chi^2(P \| Q), \chi^2(Q \| P)\}$ for two distributions $P$ and $Q$.  
In Lemmas~\ref{lm:E-J} and~\ref{lm:E-I} in Appendix~\ref{sec:proofddt}, we show how to bound the Dirichlet form of our projected chain $\overline{M}$, denoted by $\overline{\mathcal{E}}$. Compared to the approach of~\citet{ge2018simulated}, our argument for bounding  $\overline{\mathcal{E}}$ is more straightforward and yields simpler and equally tight 
bound on the Poincar\'{e} constant of the simulated tempering chain in the decomposition theorem.  
\end{remark}

We are now ready to state our main decomposition theorem, which 
relates the restricted spectral gap of $M$ to $\mathrm{SpecGap}_{\cX^0}(M_{(i, j)})$ for $i\in [L], j \in [\m]$ and $\mathrm{SpecGap}(\overline{M})$.  Roughly speaking,  this result guarantees that the restricted spectral gap of $M$ is $\Omega ( (C_1 C_2 C_3)^{-1} )$.   

\begin{theorem}\label{thm: ddt}
Consider the simulated tempering Markov chain \( M \) given in Definition \ref{def: simulated}.  
Suppose Assumptions~\ref{ass:pi},~\ref{ass:Mij} and~\ref{ass:barM} hold, and define 
\begin{equation}
    \theta = P( [L] \times \cX^0), \quad \phi = \min_{i, j} P_{(i, j)} (\cX^0).
\end{equation}
Then, we have 
$\mathrm{SpecGap}_{[L] \times \cX^0}(M) \geq C_M^{-1}$ where     
\begin{equation}\label{eq:CM}
{C_M =\max \left\{3 \theta C_3, \; \frac{\theta C_1C_2}{\phi(1 - \lambda) } \left({(2+\lambda) C_3} + 1\right)\right\}{\mathcal{E}_{[L] \times \cX^0}}(g, g)}. 
\end{equation}   
\end{theorem}
\begin{proof}
See Appendix \ref{sec:proofddt}.
\end{proof} 
{When $\cX^0 = \cX$, we obtain a bound on the spectral gap of $M,$ which is similar in spirit to the continuous-time decomposition theorem of~\citet{ge2018simulated}, though the two results are not directly comparable due the discrete-time setting and our construction of $\overline{M}$.} 

\begin{corollary}\label{coro:gap}
Consider the simulated tempering Markov chain \( M \) given in Definition \ref{def: simulated}. 
Suppose Assumptions~\ref{ass:pi},~\ref{ass:Mij} and~\ref{ass:barM} hold with $\cX^0 = \cX$. Then,  
$\mathrm{SpecGap}(M) \geq C_M^{-1}$ where  $C_M$ is given by~\eqref{eq:CM}. 
\end{corollary}
\begin{proof}
    This immediately follows from Theorem~\ref{thm: ddt}. 
\end{proof}

\subsection{Mixing Time Bounds}
Given two probability distributions $\Pi_1$ and $\Pi_2$ with densities $\pi_1$ and $\pi_2$, we define their total variation distance by $\|\Pi_1 - \Pi_2 \|_{\mathrm{tv}} = \int | \pi_1(x) - \pi_2(x) | \d x$, which takes values in $[0, 2]$.  
In addition to Assumptions~\ref{ass:pi} to~\ref{ass:barM}, we need to further require that $M$ is both reversible and lazy, where ``lazy'' means that $M( (i, x), \{ (i, x)\} ) \geq 1/2$ for any $(i, x) \in [L] \times \cX$.  Both conditions are very mild and commonly used in the sampling literature. 
Then, if $M$ has a positive spectral gap,  \(M\) converges exponentially fast to its stationary distribution \(P\) in TV distance~\citep{montenegro2006mathematical}.  
It was shown in~\citet{atchade2021approximate} that a positive restricted spectral gap can also be used to obtain an upper bound on the mixing time. The following lemma follows from the result of~\citet{atchade2021approximate}, where we also characterize the exponential convergence rate at each temperature level.

\begin{lemma}\label{cor:temp}
Assume that the simulated tempering Markov chain \( M \), defined in Definition~\ref{def: simulated}, is reversible and lazy, with stationary distribution \( P \). Let \( P^0 \) denote the initial distribution, and suppose that \( P^0 \) is absolutely continuous with respect to \( P \). Define \( f_0 \) by
\begin{equation}\label{eq:f_0}
    P^0(i, \mathrm{d}x) = f_0(i, x)\, P(i, \mathrm{d}x).
\end{equation}
  Let \(P^N\) denote the distribution of the chain after \(N\) steps, and for each temperature level \(i \in [L]\), define the conditional distribution
  $P_i^N(\d x) \propto P^N(i, \d x).$ 
  Fix $\varepsilon \in (0, 1)$. 
  Suppose that there exist constants \( B > 1 \), \(\infty \geq  q > 2 \) and a measurable set $\cX^0 \subset \cX$ such that
  \begin{enumerate} 
      \item $\| f_0 \|_{\mathcal{L}^q(P)} \leq B$,  where $\| \cdot \|_{\mathcal{L}^q(P)}$ denotes the $\mathcal{L}^q$-norm w.r.t.  $P$, 
      \item $P([L] \times \cX^0) \geq 1 -  (\frac{\varepsilon^2 }{20 B^2} )^{q/(q-2)}$,
      \item $ \mathrm{SpecGap}_{[L] \times \cX^0}(M) \geq C_M^{-1}.$
  \end{enumerate}   
Then, for  $N \geq C_M \log (2B^2 / \varepsilon^2 )$, 
we have
\begin{align}\label{eq:restricted-conv}
    \|P^N - P\|_{\mathrm{tv}} \le \varepsilon,    \text{ and }
    \|P_i^N - P_i\|_{\mathrm{tv}} \le 
    \frac{3 \varepsilon}{2 \min_{k \in [L]} r_k}   
    \text{ for all } i \in [L]. 
\end{align} 
\end{lemma}
\begin{proof}
See Appendix \ref{sec:proofddt}.
\end{proof}

 {The constant \( C_M \) controls the rate of convergence of the simulated tempering chain to its stationary distribution in total variation distance; larger values of \( C_M \) correspond to slower convergence.
}
\section{ {Analysis of the Simulated Tempering Metropolis–-Hastings Algorithm for Multivariate Gaussian Mixtures}}\label{sec: example}

Let function $f\colon \mathbb{R}^d \to \mathbb{R}$ be defined by 
\begin{equation}\label{eq:potential}
f(x) = - \log \left\{ \sum_{i=1}^{\m} w_i  e^{-\frac{1}{2} (x - \mu_i)^\top \Sigma^{-1} (x - \mu_i) } \right\}, 
\end{equation}
where $\mu_i \in \mathbb{R}^d$ for each $i$, $\Sigma \in \mathbb{R}^{d\times d}$ is a positive definite matrix, and $w_i > 0$ for each $i$ such that $\sum_{i=1}^{\m} w_i = 1$.   We want to sample from a probability distribution $P^*$ with  probability density function $p^*(x) \propto e^{-f(x)}$. The target density $p^*$ corresponds to a mixture of  $\m$ Gaussian components, where each component has mean \( \mu_i  \), weight \( w_i \), and shared covariance matrix \(\Sigma\). We assume that we do not have access to the individual means \( (\mu_i)_{i=1}^{\m} \) and weights \( (w_i)_{i=1}^{\m} \), but we can evaluate $f(x)$ at any point \( x \in \mathbb{R}^d \); the problem would become trivial if $(\mu_i, w_i)_{i=1}^{\m}$ are known. 
 
We illustrate the use of our decomposition theorem by analyzing the sampling complexity of the simulated tempering Metropolis--Hastings (STMH) algorithm for target distribution $P^*$.

\begin{definition}[STMH]\label{def:stmh}
Given a target probability density $p^*(x) \propto \exp(- f(x))$, let  
\begin{equation*}
        M^* = M\left( (p^*_i)_{i=1}^L, (r_i)_{i=1}^L, (M^*_i)_{i=1}^L, \lambda \right), 
\end{equation*} 
denote the simulated tempering Markov chain defined in Definition~\ref{def: simulated}. Here,
\begin{equation}\label{eq:pi}
    p^*_i(x) \propto \exp(-\beta_i f(x)),
\end{equation}
where \( 0 < \beta_1 <   \cdots < \beta_L = 1 \) is a sequence of inverse temperatures. The transition kernel \( M^*_i \) is that of the Metropolis--Hastings algorithm  targeting \( p^*_i \) with a symmetric Gaussian proposal density \( q(x, y) = \mathcal{N}(y; x, \eta I) \), where \( \eta > 0 \) is the step size. The constants \( r_i \) are set proportional to \( Z_i / \widehat{Z}_i \), where \( Z_i \) is the (unknown) normalizing constant of \( p^*_i \), and \( \widehat{Z}_i \) is its estimate. Specifically, we set
\begin{align}\label{eq:r_i}
r_i = \frac{Z_i / \widehat{Z}_i}{\sum_{k=1}^L Z_k / \widehat{Z}_k}
\end{align}
so that \( \sum_{i=1}^L r_i = 1 \).
\end{definition}

 {We define \( r_i \propto Z_i/\widehat{Z}_i \)  because the true normalizing constants \( Z_i \)  (also known as partition functions, where 
$Z$ is viewed as a function of the temperature index $i$)  are typically unknown in practice. When implementing the STMH algorithm, the acceptance probability in Equation~\eqref{eq:accp} is given by
\[
a\big((i, x), (i', x)\big) = \min \left\{ \frac{\widehat{Z}_i \exp(-\beta_{i'} f(x))}{\widehat{Z}_{i'} \exp(-\beta_i f(x))},\, 1 \right\}.
\]
This choice ensures that the acceptance probability depends only on the estimated normalizing constants \( \widehat{Z}_i \) and not the true values \( Z_i \), thereby making the algorithm implementable even when \( Z_i \) are  unknown.} Since acceptance probability depends only on the ratio \( \widehat{Z}_i / \widehat{Z}_{i'} \), it is sufficient to estimate the normalizing constants up to a common multiplicative factor. We set $\widehat{Z}_1 = 1$ and estimate the other normalizing constants using the inductive  strategy considered by~\citet{ge2018simulated}: for $\ell = 1, \dots, L - 1$, we run STMH using $\ell$ inverse temperatures $\beta_1 < \cdots < \beta_{\ell}$ and use the samples at temperature level $\ell$ to compute the estimate $\widehat{Z}_{\ell + 1}$. This estimation procedure is summarized in Algorithm~\ref{alg:main} in Appendix~\ref{sec:app-alg}.
After obtaining all estimates $(\widehat{Z}_i)_{i \in [L]}$, we run STMH using all $L$ temperature levels, and the samples collected at the $L$-th temperature level can be used to approximate the original target density $p^*(x) \propto e^{-f(x)}$. Algorithm~\ref{alg:simulated_tempering} in Appendix~\ref{sec:app-alg} summarizes the STMH algorithm given in  Definition~\ref{def:stmh}, assuming all partition functions are known.

The following theorem provides a mixing time bound for the STMH algorithm targeting the distribution $P^*$.  
Directly applying the spectral gap decomposition in Corollary~\ref{coro:gap} is very challenging,  since it remains unclear how to effectively bound the spectral gap of the Metropolis--Hastings chain over an unbounded domain. In contrast, such bounds are more tractable within bounded regions~\citep{yuen2000applications}. 
This motivates the use of  the decomposition of restricted spectral gap given in   Theorem~\ref{thm: ddt}.

\begin{theorem} \label{thm: pa}
Let \( f(x) \) be defined as in Equation~\eqref{eq:potential}, and define  
\[
D \coloneqq \max\left\{ \max_{k \in [\m]} \|\mu_k\|,\, \sqrt{\gamma_{\min}} \right\}, \quad
w_{\min} := \min_{1 \leq j \leq \m} w_j, \qquad 
\kappa \coloneqq \frac{\gamma_{\text{max}}}{\gamma_{\text{min}}},
\]
where \( \gamma_{\text{max}} \) and \( \gamma_{\text{min}} \) denote the largest and smallest eigenvalues of the covariance matrix \( \Sigma \), respectively. 
Then, assuming $d$ is fixed, 
STMH algorithm can be used to generate a sample from a distribution that is within \( \varepsilon \) total variation distance of the target distribution \( P^* \), in time
\[
\operatorname{poly} \left( \frac{1}{w_{\min}}, D, \kappa \right) \cdot \log^3\left(\frac{1}{\varepsilon}\right).
\] 

\end{theorem}
\begin{proof}
See Appendix \ref{sec:proof-example}, where  a precise version of the time complexity bound is also provided. 
\end{proof}

\begin{remark}
The parameter \( D \) captures the spread and separation of the Gaussian mixture components since $\max_{i \neq j} \| \mu_i - \mu_j \| \leq 2D$ by the triangle inequality. 
\end{remark}

\begin{remark}
 {The time complexity in Theorem~\ref{thm: pa} is polynomial in the separation parameter \( D \), logarithmic in the inverse accuracy \( 1/\varepsilon \), and exponential in the dimension \( d \). }
 In fixed-dimensional settings---where \( d \) is constant---this result provides high-accuracy guarantees, showing that STMH is particularly effective for sampling from such mixture distributions.   For comparison, the proximal sampler incurs a time complexity with exponential dependence on $D$, since the log-Sobolev constant and Poincar\'{e} constant for a Gaussian mixture distribution decays exponentially in \( D \); see~\citet{schlichting2019poincare}.  
To further illustrate the scaling behavior of the proximal sampler, consider a mixture of two Gaussians with identical covariance matrices and means located at  \( (-m/2,  \dots, -m / 2) \) and \( (m/2,   \dots, m/2) \), respectively.  
In this case, the separation becomes \( D = m \sqrt{d} \), which implies that the complexity of the proximal sampler  scales exponentially in both $m$ and $d$.   
 {The simulated tempering Langevin Monte Carlo (STLMC) algorithm of~\citet{ge2018simulated} admits an upper bound that scales polynomially with \(D\),  \(d\), and  \(1/\varepsilon\). In contrast, the upper bound we obtain for STMH exhibits only a logarithmic dependence on \(1/\varepsilon\), representing a state-of-the-art theoretical guarantee among existing known upper bounds.} Table~\ref{tab:comparison} summarizes the theoretical complexity of STMH in comparison with several other sampling algorithms. It can be seen that no algorithm dominates STMH in terms of the complexity dependence on $D$ and $1/\varepsilon$. 
 {Our upper bound on the time complexity of STMH has exponential dependence on $d$, which is likely due to the use of the random walk Metropolis--Hastings (RWMH) sampler for each local chain $M_i$.} 
We conjecture that by replacing it with  proximal sampler~\citep{chen2022improved,he2024a} or Metropolis-adjusted Langevin algorithm~\citep{wu2022minimax}, the resulting simulated tempering algorithm may achieve a better complexity dependence on $d$. 
\end{remark}

\begin{remark}
 {Theorem~\ref{thm: pa} establishes convergence guarantees for STMH when sampling from a mixture of Gaussians with a shared covariance matrix. This result can be naturally extended to target distributions that are sufficiently close to such mixtures, following an approach similar to that of~\citet{ge2018simulated}.
 In such cases, the time complexity would additionally depend on closeness between the actual distribution and the Gaussian mixture approximation. We also anticipate that similar techniques can be adapted to handle mixtures of log-concave distributions that are translates to each other, or more broadly, distributions that are well-approximated by such mixtures. However, as demonstrated by~\citet{ge2018simulated}, some seemingly mild violations of the assumptions, such as component covariance matrices differing by a constant factor, can lead to exponential time complexity for any reasonable algorithm with similar guarantees.}
\end{remark}
\begin{table*}[htbp]
\centering
\renewcommand{\arraystretch}{1.2}
{ 
\begin{tabular}{|c|c|c|c|}
\hline
\textbf{Algorithm} &  $d$ &  $D$ &  $1/\varepsilon$ \\ \hline
 STLMC~
\citep{ge2018simulated}  & poly & poly & poly \\ \hline
 LMC~
\citep{vempala2019rapid} & $d$ & $\exp(D^2)^2$ & $1/\varepsilon^2$ \\ \hline
Annealed LMC~\citep{guo2024provable}  & poly & poly & poly \\ \hline
Proximal Sampler~\citep{fan2023improved}  & $d^{1/2}$ & $\exp(D^2)$ & $\log(1/\varepsilon)$ \\ \hline
LMC (data-initialized)~\citep{koehler2024efficiently} & poly & $\exp(D^2)^2$ & $1/\varepsilon^2$ \\ \hline
 {STMH} &  {exp} &  {poly} & $ {\log^2(1/\varepsilon)}$  \\ \hline
\end{tabular}
} 
\caption{Dependence of time complexity on $d$, $D$ and $\varepsilon$ for sampling from densities of the form in~\eqref{eq:potential}.
}
\label{tab:comparison}
\end{table*}

The proof of Theorem \ref{thm: pa} is divided into several steps. In order to apply the decomposition arguments developed in Section~\ref{sec:prob}, we  first define an approximate STMH chain that satisfies Assumption~\ref{ass:pi}. 

\begin{definition}[Approximate STMH]\label{def:approx-stmh}
Let  
\begin{equation*}
        \widetilde{M} = M\left( (\widetilde{p}_i)_{i=1}^L, (r_i)_{i=1}^L, (\widetilde{M}_i)_{i=1}^L, \lambda \right), 
\end{equation*} 
denote the simulated tempering Markov chain defined in Definition~\ref{def: simulated}, where 
\begin{equation}\label{eq:tilde-pi}
    \widetilde{p}_i(x) \propto \sum_{j=1}^\m w_j \exp\left\{ -\frac{\beta_i}{2} (x - \mu_j)^\top \Sigma^{-1} (x - \mu_j) \right\},
\end{equation}
and transition kernel \( \widetilde{M}_i \) is that of the Metropolis--Hastings algorithm targeting \( \widetilde{p}_i \), with a {symmetric Gaussian proposal density} \( q(x, y) = \mathcal{N}(y; x, \eta I) \) where step size \( \eta > 0 \). The weights \( (r_i)_{i = 1}^L \) and $\lambda$ are the same as  in Definition~\ref{def:stmh}.
\end{definition}

The stationary densities $ \widetilde{p}_i$ in Definition~\ref{def:approx-stmh} are mixtures of Gaussian densities; denote the component distributions by $\widetilde{P}_{(i, j)}$. 
This enables us to apply Theorem~\ref{thm: ddt}. It is important to note that this approximate STMH sampler cannot be implemented in practice, as it requires knowledge of the component distributions, which is typically unavailable.  
We show that, for some $\cX^0 \subset \mathbb{R}^d$, the $\cX^0$-restricted spectral gap of the STMH chain constructed in Definition~\ref{def:stmh} is comparable to that of the approximate STMH chain given in Definition~\ref{def:approx-stmh} (see Lemma~\ref{lm: var} in Appendix~\ref{sec:appx-proof-compare}). This comparison result is a discrete-time extension of the argument used in \citet{ge2018simulated}, which suggests that adjusting the temperature is roughly equivalent to modifying the variance of a Gaussian distribution. 
As a result, it suffices to obtain a lower bound on the $\cX^0$-restricted spectral gap for the approximate STMH chain.  Note that since $\beta_L = 1$,  \( p_L = \widetilde{p}_L = p^* \).

Next, we apply Theorem~\ref{thm: ddt} to derive a lower bound on the \( \cX^0 \)-restricted spectral gap of the approximate STMH chain. This involves three  steps.  First, we verify inequality~\eqref{eq:decomp-c1} in Appendix~\ref{sec:appx-proof-c1} by comparing the transition density of the chain $\widetilde{M}_{(i, j)}$ with that of $\widetilde{M}_i$.  
Second, we need to compute the $\cX^0$-restricted spectral gap of the ``local chain'' $\widetilde{M}_{(i, j)}$, which we define as the random walk Metropolis--Hastings algorithm targeting $\widetilde{P}_{(i, j)}$, the Gaussian distribution with mean $\mu_j$ and covariance matrix $\beta_i^{-1} \Sigma$. It is well known that for strongly log-concave target distributions, the Metropolis--Hastings algorithm should mix fast~\citep{johndrow2018fast}. To explicitly compute the bound, we apply the path method on continuous spaces proposed by~\citep{yuen2000applications}, which is detailed  in Appendix~\ref{sec:appx-proof-c2}. 
Finally, the projected chain captures transitions between mixture components and temperature levels. Intuitively, it should mix fast because (i) if $\beta_1$ is sufficiently small, the component distributions  $\widetilde{P}_{(1, 1)}, \dots, \widetilde{P}_{(1, \m)}$ overlap significantly, and (ii) if $\beta_i/\beta_{i-1}$ is not too big, then $\widetilde{P}_{(i, j)}$ and $\widetilde{P}_{(i-1, j)}$ also share substantial overlap. To compute a lower bound on the spectral gap  for the projected chain, we apply the well-known canonical path method~\citep{levin2017markov}; see   Appendix~\ref{sec:appx-proof-c3}.   

To conclude the proof, we derive error bounds on the estimated partition functions in Appendix~\ref{sec:main_thm}.

\section{Simulation Study}\label{sec:simulation}
To numerically investigate the complexity of the STMH algorithm, we perform a simulation study with target distribution being a symmetric two-dimensional Gaussian mixture distribution, whose density is given by
\[
p^*(x) = \frac{1}{2} \, \mathcal{N}\left(x; -\frac{D}{2\sqrt{2}} \cdot {1_2}, I_2\right) + \frac{1}{2} \, \mathcal{N}\left(x; \frac{D}{2\sqrt{2}} \cdot {1_2}, I_2\right),
\]
where \(I_2\) is the \(2 \times 2\) identity matrix and \({1_2} = (1, 1)^\top\). The parameter \(D\) controls the separation between the two components. We vary the parameter \(D\) to explore how increasing the separation between the modes affects the convergence behavior of the algorithm. 
For each value of \(D\), we run the STMH algorithm with the parameters specified in Appendix~\ref{sec:ass_para} and initialized at $(10, 10)$. 
To assess convergence, we monitor how quickly the empirical mean of the samples, denoted by $\hat{\mu}$, approaches the true mean of the target distribution, \((0, 0)\). If the chain has not yet mixed, it tends to spend more time near one mode, resulting in an empirical mean that deviates from the target mean.
Moreover, according to ~\citet{nishiyama2022lower},  when the empirical and target distribution have the same covariance matrix, 
we can lower bound their TV distance by \( \|\hat{\mu}\|^2 / (C + \|\hat{\mu}\|^2) \) for some constant $C> 0$, 
which is a monotone increasing function of \( \|\hat{\mu}\|^2 \). This justifies the use of \( \|\hat{\mu}\| \) as a measure of convergence.  
To reduce variability, we repeat the simulation $20,000$ times for each fixed $D$ and average the empirical means over all runs. {To provide a benchmark, we also run the baseline Metropolis–Hastings (MH) algorithm under the same setup and compare its convergence behavior with that of STMH.} All simulations were performed on a standard consumer-grade CPU with parallelization and completed within approximately six hours.

\begin{figure}[ht]     
    \centering   \includegraphics[width=0.48\linewidth]{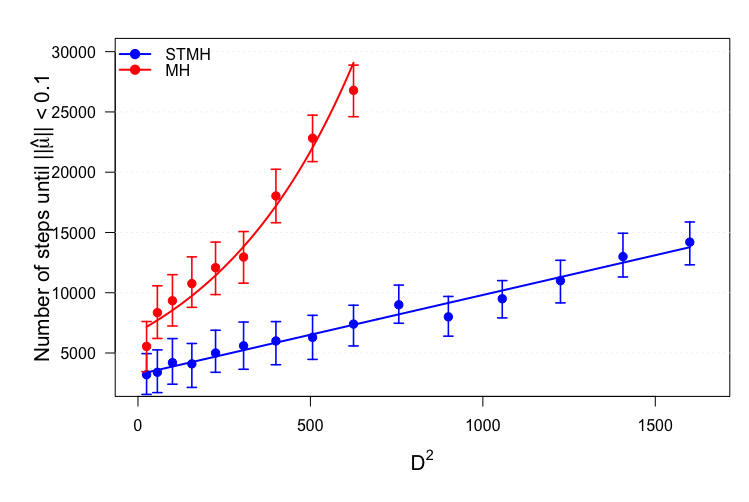}    \includegraphics[width=0.48\linewidth]{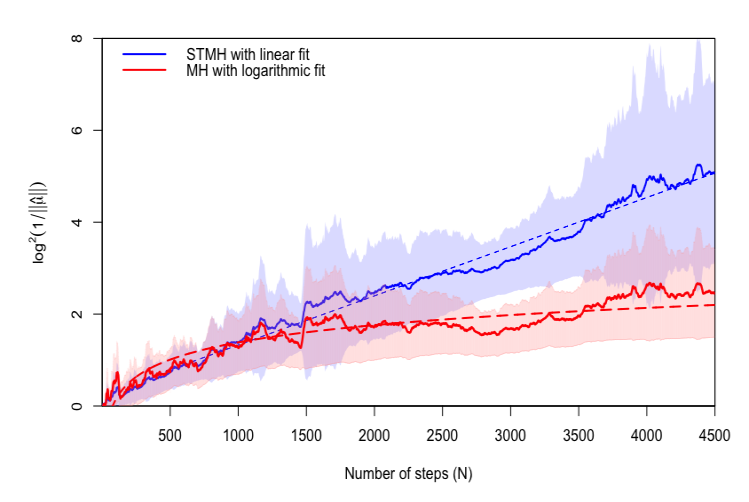}
    \captionof{figure}{Left: number of steps until \( \|\widehat{\mu}\| < 0.1 \) versus   \( D^2 \). Right: \( \log^2\left(1 / \|\widehat{\mu}\|\right) \) versus the number of steps \( N \) for \( D = 30 \).} \label{fig:combined}
\end{figure}

In the left panel of Figure~\ref{fig:combined}, we plot the number of steps required for the empirical mean to fall below 0.1 for both STMH and the baseline MH algorithm. The error bars are obtained by computing the 95\% coverage interval of the empirical mean at each step, and then determining the corresponding number of steps needed for the lower and upper bounds of the interval to cross the 0.1 threshold. The results show that this number grows approximately linearly with \(D^2\) for STMH, consistent with the theoretical bound in Theorem~\ref{thm: pa}. In contrast, for the baseline MH algorithm the number of steps increases exponentially with \(D^2\), reflecting its slower convergence, which is consistent with the well-known behavior of MH. To analyze how the mixing time depends on the threshold \(\varepsilon\), we fix \(D = 30\) and plot how $\|\hat{\mu}\|$ varies with the number of steps \(N\) in the right panel of Figure~\ref{fig:combined}. For STMH, the results  show an approximately linear relationship between the number of steps and $\log^2\!\big(1/\|\hat{\mu}\|\big)$, consistent with the theoretical bound in Theorem~\ref{thm: pa}. In contrast, for the baseline MH algorithm the relationship is closer to logarithmic, implying a polynomial dependence of the mixing time on $1/\varepsilon$, which aligns with the well-known limitations of MH in achieving high accuracy. The error bars represent 95\% coverage intervals of the empirical mean, with the delta method applied to obtain the corresponding intervals for $\log^2\!\big(1/\|\hat{\mu}\|\big)$.  Overall, these simulations demonstrate that STMH provides more efficient sampling from this Gaussian mixture distribution compared to MH.

\section{Concluding Remarks}\label{sec:conclusion}

Simulated tempering addresses the challenge of sampling from multimodal distributions. In this work, we develop a general theoretical framework for analyzing simulated tempering and demonstrate its effectiveness through a detailed analysis of simulated tempering combined with the Metropolis–Hastings algorithm for sampling from Gaussian mixtures. Our framework can be used to analyze simulated tempering combined with other local MCMC samplers, such as the Metropolis-adjusted Langevin algorithm (MALA) and proximal algorithms, but verifying  Assumptions~\ref{ass:Mij} and~\ref{ass:barM} and computing the constants $C_1, C_2, C_3$ in these settings may be more involved.
In particular, since it has been shown in the literature that the dimensional dependence of the complexity of MALA for log-concave target distributions is $\tilde{\Theta}(\sqrt{d})$~\citep{chewi2021optimal, wu2022minimax}, it would be interesting to investigate if the complexity of simulated tempering combined with MALA has a similar polynomial dependence on $d$,  improving on the exponential dependence in our result.  
Another  promising direction for future work is to generalize our techniques to a broader class of target distributions beyond Gaussian mixtures. The argument of~\citet{ge2018simulated} can be used to extend our result to target distributions that are ``sufficiently close'' to Gaussian mixtures. 
Finally, 
while this work primarily focuses on establishing lower bounds on the spectral gap, an interesting direction for future work is to investigate the tightness of these bounds by also deriving upper bounds on the spectral gap for simulated tempering Markov chains. 

\section{Acknowledgments}
Jhanvi Garg and Quan Zhou were supported in part by NSF grants DMS-2245591 and DMS-2311307. 
Krishnakumar Balasubramanian was supported in part by NSF grant DMS-2413426. 
\bibliographystyle{plainnat}

\clearpage 
\newpage 

\section*{\large{Appendices}}

\appendix
\addcontentsline{toc}{section}{Supplementary Material}
\renewcommand{\thesection}{\Alph{section}}

\section{Simulated Tempering Metropolis--Hastings algorithm}\label{sec:app-alg}

\begin{algorithm}
\caption{Simulated Tempering Metropolis--Hastings Algorithm}
\label{alg:simulated_tempering}
\begin{algorithmic}[1]
\State \textbf{Input:} function $f$, inverse temperatures $\beta_1, \ldots, \beta_{\ell}$, partition function estimates $\widehat{Z}_1, \ldots, \widehat{Z}_{\ell}$, number of steps $N$, step size $\eta$, rate $\lambda$, initial covariance matrix $\Sigma_0$.
\State Sample $x_0 \sim \mathcal{N}(0, \Sigma_0)$ 
\State $i \gets 1, x \gets x_0, n \gets 0$ 
\While{$n < N$}  
    \State Sample $u \sim \text{Bernoulli}(\lambda)$.
    \If{$u = 0$} 
\State Propose \( x' \sim \mathcal{N}(x, \eta I) \)
\State Sample \( v \sim \mathrm{Uniform}(0,1) \)
        \If{ \( v < \min\left\{1, \frac{e^{-\beta_i f(x')} }{e^{-\beta_i f(x)} }\right\} \) }
            \State  \( x \gets x' \)
        \EndIf
    \Else
        \State  Propose $i' = i \pm 1$, each with probability $1/2$ 
        \If{$1 \leq i' \leq \ell$}
            \State Sample $v \sim \mathrm{Uniform}(0, 1)$
            \If{ $v <  \min \left\{1,  \frac{e^{-\beta_{i'} f(x)} / \widehat{Z}_{i'}}{e^{-\beta_i f(x)} / \widehat{Z}_i} \right\}$} 
                \State $i \gets i'$ 
            \EndIf
        \EndIf
    \EndIf
    \State  $n \gets n + 1$ 
\EndWhile
\State  {\textbf{Output:} Sample \((x, i)\) collected at the \(N^{\text{th}}\) step.}
\end{algorithmic}
\end{algorithm}

\begin{algorithm}[H]
\caption{Partition Function Estimation}
\label{alg:main}
\begin{algorithmic}
\State  \textbf{Input:} function $f$, inverse temperature sequence  $\beta_1 < \cdots < \beta_{L}$ and  {number of samples $s$}.
\State  $\widehat{Z}_1 \gets 1$ 
\For{$\ell = 1$ \textbf{to} $L$}
\State  {Repeat Algorithm~\ref{alg:simulated_tempering} until $s$ samples $(x_j)_{j=1}^s$ are obtained at temperature level $\ell$}.
\State  $\widehat{Z_{\ell +1}} \gets  ( \widehat{Z_{\ell }} / s )  \sum_{j=1}^s e^{\left(-\beta_{\ell +1}+\beta_{\ell }\right)f(x_j)} $ 
\EndFor 
\end{algorithmic}
\end{algorithm}
\begin{remark}
 {Algorithm~\ref{alg:simulated_tempering} is always run for a fixed number of steps \( N \) and returns the sample obtained at the final step. In Algorithm~\ref{alg:main}, if this sample is not from the desired temperature level, Algorithm~\ref{alg:simulated_tempering} is simply re-run for another \( N \) steps.}
\end{remark}
\newpage

\section{Proofs for Section~\ref{sec:prob}} \label{sec:proofddt}

\subsection{Proof of  Lemma~\ref{lm:barM}}
 
Clearly, $\sum_{i, j} \overline{P}(i, j) = 1$. Hence, it only remains to   check the detailed balance condition 
\[
\overline{P}((i, j)) \, \overline{M}((i, j), (i', j')) = \overline{P}((i', j')) \, \overline{M}((i', j'), (i, j)) 
\] 
for two types of moves. First, let $i$ be fixed and consider $j \neq j'$. Then, 
\begin{align*}
 r_i w_{(i, j)} P_{(i, j)}(\cX^0) \overline{M}( (i, j), (i, j')) 
&= r_i (1-\lambda)   \int_{\cX^0}  w_{(i, j)} p_{(i,j)}(x)
\frac{ w_{(i, j')} p_{(i, j')}(x) }{p_i(x)} \d x, 
\end{align*}
which is clearly symmetric with respect to $j$ and $j'$. 
Similarly, if $j$ is fixed and $i' = i \pm 1$ (assuming both $i, i' \in [L]$), we have 
\begin{align*}
  r_i w_{(i, j)} P_{(i, j)}(\cX^0) \overline{M}( (i, j), (i', j))     
  &= \frac{\lambda}{2}\int_{\cX^0}  r_i w_{(i, j)} p_{(i,j)}(x) a( (i, j, x), (i', j, x) )  \d x \\
  &= \frac{\lambda}{2}\int_{\cX^0}  \min \left\{ r_i w_{(i, j)} p_{(i,j)}(x) , r_{i'} w_{(i', j)} p_{(i', j)}(x)   \right\}  \d x, 
\end{align*}
which is symmetric with respect to $i$ and $i'$.

\subsection{Proof of Theorem~\ref{thm: ddt}}

We first prove an auxiliary lemma on the Dirichlet form of the simulated tempering Markov chain. 

\begin{lemma} 
\label{lm: dirform}
The \([L] \times \cX^0\)-restricted Dirichlet form \(\mathcal{E}_{[L] \times \cX^0}\) of the simulated tempering Markov chain \(M\), defined in Definition~\ref{def: simulated}, can be expressed by 
\begin{align*}
  \mathcal{E}_{[L] \times \cX^0}(g, g) 
  = (1 - \lambda) \sum_{i=1}^L r_i\, \mathcal{E}_{i, \cX^0}(g_i, g_i) 
  + \lambda\, \mathcal{E}^I_{\cX^0}(g, g),
\end{align*}
where \(\mathcal{E}_{i,\cX^0}\) is the \(\cX^0\)-restricted Dirichlet form of the Markov chain \(M_i\), \(g \in \mathcal{L}^2([L] \times \cX, P)\) with \(g_i(x) = g(i, x)\) for each \(i \in [L]\) and
\begin{align*}
  \mathcal{E}^I_{\cX^0}(g, g) 
  = \frac{1}{4} \sum_{\substack{i, i' \in [L] \colon i' = i \pm 1}} 
  \int_{\cX^0} \bigl(g(i, x) - g(i', x)\bigr)^2 
  r_i\, p_i(x)\, a\bigl((i, x), (i', x)\bigr)\, \mathrm{d}x.
\end{align*}
\end{lemma}

\begin{proof}
Since the stationary density of $M$ is $p(i, x) = r_i p_i(x)$ and either $x$ or $i$ is fixed in each simulated tempering iteration, the restricted Dirichlet form \(\mathcal{E}_{[L] \times \cX^0}(g, g)\) can be expressed by 
\begin{equation}\label{eq:dirichlet}
\begin{aligned}
\mathcal{E}_{[L] \times \cX^0}(g, g) &=   
\frac{1}{2} \, \sum_{i=1}^{L}
\int_{\cX^0} \!\!\int_{\cX^0}
\bigl(g(i, x) - g(i, y)\bigr)^2 \,
r_i \, p_i(x)\, 
M\bigl((i, x), (i, \d y)\bigr)\,
\mathrm{d}x  \\
&\quad +  \frac{1}{2}   \sum_{i, i' \in [L] \colon i' = i \pm 1 }
\int_{\cX^0}
\bigl(g(i, x) - g(i', x)\bigr)^2 \,
r_i \, p_i(x)\, 
M\bigl((i, x), (i', x)\bigr)\,
\mathrm{d}x.  
\end{aligned}
\end{equation}
Since $ M\bigl((i, x), (i, \d y)\bigr)  = (1 - \lambda) M_i (x, \d y)$ and $M_i$ has stationary density $p_i$, 
\begin{align*}
&    \frac{1}{2} \, 
\sum_{i=1}^{L}
\int_{\cX^0}\!\!\int_{\cX^0}
\bigl(g(i, x) - g(i, y)\bigr)^2 \,
r_i \, p_i(x)\, 
M\bigl((i, x), (i, \d y)\bigr)\,
\mathrm{d}x   \\
&= \frac{1 - \lambda}{2}
\sum_{i=1}^{L}
r_i \int_{\cX^0}\!\!\int_{\cX^0}
\bigl(g_i(x) - g_i(y)\bigr)^2 \,
 p_i(x)\, 
M_i\bigl(  x, \d y  \bigr)  \d x = (1 - \lambda) \sum_{i=1}^{L} r_i \,\mathcal{E}_{i, \cX^0} \bigl(g_i, g_i\bigr). 
\end{align*}
For the second term on the right-hand side of~\eqref{eq:dirichlet}, note that 
\begin{align*}
   M\bigl((i, x), (i', x)\bigr) = \frac{\lambda}{2} a((i, x), (i', x)),
\end{align*}
where the acceptance probability function $a$ is given by~\eqref{eq:accp}. 
A straightforward calculation then concludes the proof of the lemma.
\end{proof}

Next, we prove two key lemmas about the Dirichlet form $\overline{\mathcal{E}}$ of the Markov chain $\overline{M}$ constructed in Definition~\ref{def:project-chain}. Let $\theta = P( [L] \times \cX^0)$. 
Recall that under Assumption~\ref{ass:pi}, we can augment the stationary density to 
\begin{equation*} 
    p(i, j, x) =  r_i w_{(i,j)} p_{(i,j)}(x). 
\end{equation*} 
We still denote the corresponding probability measure by $P$. 
Let $P_0$ denote the conditional probability measure given  $X \in \cX^0$, whose density is given by 
\begin{equation}
    p_0(i, j, x) =  \frac{ r_i w_{(i,j)} p_{(i,j)}(x) }{ \theta }. 
\end{equation}
The Dirichlet form $\overline{\mathcal{E}}$  can be expressed by 
\begin{equation}\label{eq:decompose-bar-E}
    \overline{\mathcal{E}}(\overline{g}, \overline{g}) =   \overline{\mathcal{E}}^J(\overline{g}, \overline{g}) + 
    \overline{\mathcal{E}}^I(\overline{g}, \overline{g}), 
\end{equation}
where 
\begin{align*}
\overline{\mathcal{E}}^J(\overline{g}, \overline{g}) &= \frac{1}{2 \theta} \sum_{i=1}^{L}    \sum_{j, j'=1}^{\m} (\overline{g}(i,j) - \overline{g}(i,j'))^2 \, r_iw_{(i,j)} P_{(i, j)} (\cX^0) \, \overline{M}((i,j), (i,j')), \\ 
\overline{\mathcal{E}}^I (\overline{g}, \overline{g}) &= \frac{1}{2 \theta}   \sum_{j=1}^{\m} \sum_{i, i' \in [L] \colon i' = i \pm 1 } (\overline{g}(i,j) - \overline{g}(i',j))^2 \, r_i w_{(i,j)} P_{(i, j)} (\cX^0) \, \overline{M}((i,j), (i',j)). 
\end{align*}

\begin{lemma}\label{lm:E-J}
Suppose Assumption~\ref{ass:pi} holds. For any $g \in \mathcal{L}^2( [L] \times \cX, P)$, define $g_i \colon \cX \rightarrow \mathbb{R}$ by $g_i(x) = g(i, x)$, and define 
\(\overline{g} \colon [L] \times [\m] \rightarrow \mathbb{R} \) by
\begin{equation*}
    \overline{g}(i,j) = \int_{\cX^0} g(i, x) \frac{ p_{(i, j)}(x) }{P_{(i, j)}(\cX^0)}   \d x. 
\end{equation*}  
Then, 
\begin{equation}
\overline{\mathcal{E}}^J(\overline{g}, \overline{g})  \leq  \frac{ 2(1 - \lambda) }{\theta}
 \sum_{i = 1}^L \sum_{j=1}^{\m}    \frac{ r_i w_{(i,j)} }{   P_{(i, j)}(\cX^0)  }   \; \mathrm{Var}_{ P_{(i, j)}, \cX^0} ( g_i ).   
\end{equation}
\end{lemma}
\begin{proof}
For every \( x \in \cX \), \(i \in [L]\) and every pair \( j, j' \in [\m] \), the following inequality holds
\begin{equation*}\label{eq:dirichlet-E-J}
  \bigl( \overline{g}(i,j) - \overline{g}(i,j') \bigr)^2 
\leq 2 \left[
\bigl(  \overline{g}(i,j) - g(i, x) \bigr)^2 
+ \bigl( \overline{g}(i,j') - g(i, x) \bigr)^2
\right].  
\end{equation*} 
Hence, using the expression for $\overline{M}( (i, j), (i, j'))$, we get 
\begin{align}
\overline{\mathcal{E}}^J(\overline{g}, \overline{g}) 
&= \frac{1 - \lambda}{2 \theta}\sum_{i = 1}^L  \sum_{j, j'=1}^{\m} (\overline{g}(i,j) - \overline{g}(i,j'))^2 r_i w_{(i,j)} \int_{\cX^0}  p_{(i,j)}(x)  p_i(j' | x) \, \d x  \\ 
&\leq  (1 - \lambda)\sum_{i = 1}^L  \sum_{j, j'=1}^{\m} r_i w_{(i,j)} \int_{\cX^0} 
\left[
\bigl( \overline{g}(i, j) - g(i, x) \bigr)^2 
+ \bigl( \overline{g}(i, j') - g(i, x) \bigr)^2 
\right]  \frac{ p_{(i,j)}(x) }{\theta}  
 p_i(j' | x) \, \d x  \\
 &= (1 - \lambda) \E_{\tilde{P}} \left[
\bigl( \overline{g}(I, J) - g(I, X) \bigr)^2 
+ \bigl( \overline{g}(I,  J') - g(I, X) \bigr)^2 
\right],  \label{eq:E-J-1}
\end{align}
where $\tilde{P}$ denotes the joint probability measure of $(I, J, J', X)$ with density 
\begin{align*}
    \tilde{p}(i, j, j', x) = \frac{ r_i w_{(i,j)} p_{(i,j)}(x) }{ \theta } p_i(j' | x) = \frac{ r_i w_{(i,j)} w_{(i, j')} p_{(i,j)}(x) p_{(i,j')}(x) }{ \theta  \, p_i(x) }, 
\end{align*}
for $i \in [L], j, j' \in [\m], x \in \cX^0$. 
Hence, under $\tilde{P}$,  the joint distribution of $(I, J, X)$ and that of $(I, J', X)$ are both given by $P_0$, and thus 
\begin{equation} \label{eq:E-J-2}
\E_{\tilde{P}} \left[
\bigl( \overline{g}(I, J) - g(I, X) \bigr)^2 
+ \bigl( \overline{g}(I,  J') - g(I, X) \bigr)^2 
\right]
= 2 \E_{P_0} \left[
\bigl( \overline{g}(I, J) - g(I, X) \bigr)^2 \right]. 
\end{equation}
Since $\overline{g}(i, j) = \E_{P_0}[  g(I, X) \mid I = i, J = j]$, we find that 
\begin{align}
   & \E_{P_0} \left[ \bigl( \overline{g}(I, J) - g(I, X) \bigr)^2 \right]
 = \E_{P_0} \left[ \mathrm{Var}_{P_0} ( g(I, X) \mid I, J ) \right] \\
&= \sum_{i = 1}^L \sum_{j=1}^{\m}  \frac{ r_i w_{(i,j)} P_{(i, j)}(\cX^0) } {\theta} \; \mathrm{Var}_{P_0} ( g(I, X) \mid I = i, J =j) \\
& = \sum_{i = 1}^L \sum_{j=1}^{\m}    \frac{ r_i w_{(i,j)} }{ \theta \, P_{(i, j)}(\cX^0)  }   \; \mathrm{Var}_{ P_{(i, j)}, \cX^0} ( g_i ),   \label{eq:E-J-3}
\end{align}
where in the last step we have used 
\begin{align*}
    \mathrm{Var}_{P_0} ( g(I, X) \mid I = i, J =j)  
= \frac{1}{2} \int_{\cX^0 \times \cX^0} [ g(i, x) - g(i, y) ]^2 
\frac{ p_{(i, j)}(x) }{  P_{(i, j)}(\cX^0) } \frac{ p_{(i, j)}(y) }{  P_{(i, j)}(\cX^0) } \d x \d y. 
\end{align*}
The claim then follows from~\eqref{eq:E-J-1},~\eqref{eq:E-J-2} and~\eqref{eq:E-J-3}. 
\end{proof}

\begin{lemma}\label{lm:E-I}
Consider the setting of Lemma~\ref{lm:E-J}. We also have 
\begin{equation}
        \overline{\mathcal{E}}^I(\overline{g}, \overline{g})  
     {\leq \frac{3 \lambda}{\theta}  \sum_{i = 1}^L \sum_{j=1}^{\m}    \frac{ r_i w_{(i,j)} }{   P_{(i, j)}(\cX^0)  }   \; \mathrm{Var}_{ P_{(i, j)}, \cX^0} ( g_i )
    + \frac{3 \lambda}{\theta} \mathcal{E}^I_{\cX^0}(g, g)}.  
\end{equation}
\end{lemma}

\begin{proof} 
For every \( x \in \cX \), every pair \( i, i' \in [L] \), and each \( j \in [\m] \),  
\[
\bigl( \overline{g}(i, j) - \overline{g}(i', j) \bigr)^2 
\leq 3 \Bigl[
  \bigl( \overline{g}(i, j) - g(i, x) \bigr)^2 
  + \bigl( g(i, x) - g(i', x) \bigr)^2 
  + \bigl( \overline{g}(i', j) - g(i', x) \bigr)^2
\Bigr].
\] 
Then, using the definition of $\overline{M}((i, j), (i', j))$ for $i' = i \pm 1$, we get  
\begin{align}
    \overline{\mathcal{E}}^I(\overline{g}, \overline{g}) 
& = \frac{  \lambda}{4 \theta}  \sum_{i, i' \in [L] \colon i' = i \pm 1 }\sum_{j=1}^{\m}  (\overline{g}(i,j) - \overline{g}(i',j))^2 r_i w_{(i,j)} \int_{\cX^0} p_{(i, j)}(x) a( (i, j, x), (i', j, x)) \d x \\
&  {\leq \frac{ 3 \lambda}{2 } \E_{\tilde{P}}  \Bigl[
  \bigl( \overline{g}(I, J) - g(I, X) \bigr)^2 
  + \bigl( g(I, X) - g(I', X) \bigr)^2 
  + \bigl( \overline{g}(I', J) - g(I', X) \bigr)^2
\Bigr]},   
\end{align} 
where $\tilde{P}$ is the probability measure of $(I, I', J, X)$ with density 
\begin{equation}
\tilde{p}(i, i', j, x) = 
\begin{cases}
\displaystyle \frac{r_i w_{(i, j)} p_{(i, j)}(x)\, a\big( (i, j, x), (i', j, x) \big)}{2\theta}, & \text{if } i' = i \pm 1, \\[8pt]
1 - \tilde{p}(i, i+1, j, x) - \tilde{p}(i, i-1, j, x), & \text{if } i' = i, \\[4pt]
0, & \text{otherwise.}
\end{cases}
\end{equation}
That is, we first draw $I, J, X  \sim P_0$ and then update $I'$ by proposing $I' = I \pm 1$ with equal probability and accept it with probability $ a( (i, j, x), (i', j, x))$. Since  the update for $I'$ given $I, J, X$ is reversible with respect to $P_0$, we also have  $(I', J, X) \sim P_0$.  
Hence, 
\begin{align*}
   &\E_{\tilde{P}}  \Bigl[
  \bigl( \overline{g}(I, J) - g(I, X) \bigr)^2  \Bigr]   =    \E_{\tilde{P}}  \Bigl[ \bigl( \overline{g}(I', J) - g(I', X) \bigr)^2  \Bigr]  
   =    \E_{P_0} \left[ \bigl( \overline{g}(I, J) - g(I, X) \bigr)^2 \right] 
\end{align*} 
which has been characterized in~\eqref{eq:E-J-3}. 
Finally, 
\begin{align*}
&\E_{\tilde{P}}  \Bigl[  \bigl( g(I, X) - g(I', X) \bigr)^2 \Bigr] \\
&= \frac{1}{2\theta}  \sum_{i, i' \in [L] \colon i' = i \pm 1 } \sum_{j=1}^\m \int_{\cX^0} (g(i, x) - g(i', x) )^2 r_i w_{(i, j)}   p_{(i, j)}(x)    a( (i, j, x), (i', j, x) )  \d x \\
& = \frac{1}{ 2\theta}  \sum_{i, i' \in [L] \colon i' = i \pm 1 } \sum_{j=1}^\m \int_{\cX^0} (g(i, x) - g(i', x) )^2 \min\left\{ r_i w_{(i,j)} p_{(i,j)}(x),\, r_{i'} w_{(i',j)} p_{(i',j)}(x) \right\}   \d x \\ 
& \leq \frac{1}{2\theta}  \sum_{i, i' \in [L] \colon i' = i \pm 1 }  \int_{\cX^0} (g(i, x) - g(i', x) )^2 \min\left\{ r_i   p_{i}(x),\, r_{i'}  p_{i'}(x) \right\}   \d x \\ 
& = \frac{2}{\theta} \mathcal{E}^I_{\cX^0}(g, g), 
\end{align*}
where $\mathcal{E}^I_{\cX^0}(g, g)$ is defined in Lemma~\ref{lm: dirform}. Note that in the inequality above, we have used  that $\sum_j \min\{ a_j, b_j \}  \leq \min \{ \sum_j a_j, \sum_j  b_j\}$ for two non-negative sequences $a_j, b_j$. 
\end{proof}

\begin{proof}[Proof of Theorem~\ref{thm: ddt}]
Fix an arbitrary $g \in \mathcal{L}^2( [L] \times \cX, P)$. 
Define, for each $i$,  $g_i \colon \cX \rightarrow \mathbb{R}$ by $g_i(x) = g(i, x)$, and 
\(\overline{g} \colon [L] \times [\m] \rightarrow \mathbb{R} \) by
\begin{equation*}
    \overline{g}(i,j) = \int_{\cX^0} g(i, x) \frac{ p_{(i, j)}(x) }{P_{(i, j) (\cX^0)}} \d x. 
\end{equation*}   
Note that $\overline{g}(i,j)$ is the conditional expectation of $g(I, X)$ given $I = i$ and $J = j$ under the joint probability measure $P_0$, 
and $\overline{P}(i, j)$ is the marginal probability of $I = i, J = j$ under $P_0$. 
Hence, by the law of total variance,  Assumption~\ref{ass:barM} and Equation~\ref{eq:E-J-3}, we find that 
\begin{align}
  \operatorname{Var}_{P_0}(g) &=  
   \operatorname{Var}_{\overline{P}}(\overline{g}) +   \sum_{i=1}^{L} \sum_{j=1}^{\m} \overline{P}(i, j) \operatorname{Var}_{P_0}(g(I, X) \mid I = i, J = j)  \\
   & \leq C_3 \overline{\mathcal{E}}(\overline{g}, \overline{g}) +
   \sum_{i = 1}^L \sum_{j=1}^{\m}    \frac{ r_i w_{(i,j)} }{ \theta  P_{(i, j)}(\cX^0)  }   \; \mathrm{Var}_{ P_{(i, j)}, \cX^0} ( g_i ). 
\end{align} 
Using $ P_{(i, j)}(\cX^0) \geq \phi$  and Assumption~\ref{ass:Mij}, 
\begin{equation}
     \sum_{i = 1}^L \sum_{j=1}^{\m}    \frac{ r_i w_{(i,j)} }{  \theta  P_{(i, j)}(\cX^0)  }   \; \mathrm{Var}_{ P_{(i, j)}, \cX^0} ( g_i ) 
    \leq \frac{C_2 }{\theta \phi}  \sum_{i = 1}^L r_i \sum_{j=1}^{\m}    w_{(i,j)}   
    \mathcal{E}_{(i, j), \cX^0}(g_i, g_i) 
    \leq \frac{C_1 C_2 }{\theta \phi}  \sum_{i = 1}^L r_i  \mathcal{E}_{i, \cX^0}(g_i, g_i). 
\end{equation}
Recall that $    \overline{\mathcal{E}}(\overline{g}, \overline{g}) =   \overline{\mathcal{E}}^J(\overline{g}, \overline{g}) + 
    \overline{\mathcal{E}}^I(\overline{g}, \overline{g}).$ By Lemma~\ref{lm:E-J}, 
\begin{equation}
\overline{\mathcal{E}}^J(\overline{g}, \overline{g})  \leq  \frac{ 2(1 - \lambda) }{\theta}
 \sum_{i = 1}^L \sum_{j=1}^{\m}    \frac{ r_i w_{(i,j)} }{   P_{(i, j)}(\cX^0)  }   \; \mathrm{Var}_{ P_{(i, j)}, \cX^0} ( g_i )
 \leq \frac{ 2(1 - \lambda) C_1 C_2 }{\theta \phi } \sum_{i = 1}^L r_i 
 \mathcal{E}_{i, \cX^0}(g_i, g_i).  
\end{equation}
By Lemma~\ref{lm:E-I}, 
\begin{align*}
        \overline{\mathcal{E}}^I(\overline{g}, \overline{g})  
    & {\leq \frac{3 \lambda}{\theta}  \sum_{i = 1}^L \sum_{j=1}^{\m}    \frac{ r_i w_{(i,j)} }{   P_{(i, j)}(\cX^0)  }   \; \mathrm{Var}_{ P_{(i, j)}, \cX^0} ( g_i )
    + \frac{3 \lambda}{\theta} \mathcal{E}^I_{\cX^0}(g, g)} \\  
    & {\leq  \frac{3 \lambda C_1 C_2 }{\theta \phi }  \sum_{i = 1}^L r_i 
 \mathcal{E}_{i, \cX^0}(g_i, g_i)  
    + \frac{3 \lambda}{ \theta} \mathcal{E}^I_{\cX^0}(g, g)} 
\end{align*}
Hence,
\begin{align*}
      &  {\frac{1}{\theta^2} \operatorname{Var}_{P, [L] \times \cX^0}(g)   = \operatorname{Var}_{P_0}(g)  
        \leq \frac{3 \lambda C_3 }{ \theta} \mathcal{E}^I_{\cX^0}(g, g) 
       + \frac{C_1 C_2  \left[ (2+\lambda)C_3 + 1 \right]  }{\theta \phi }  \sum_{i = 1}^L r_i 
 \mathcal{E}_{i, \cX^0}(g_i, g_i)}. 
\end{align*} 
Comparing with Lemma \ref{lm: dirform}, we obtain the Poincar\'{e} inequality for $M$
\begin{align*}
  \operatorname{Var}_{P, [L] \times \cX^0}(g) &  {\leq \max \left\{3 \theta C_3, \; \frac{\theta C_1C_2}{\phi(1 - \lambda) } \left({(2+\lambda) C_3} + 1\right)\right\}{\mathcal{E}_{[L] \times \cX^0}}(g, g)}
\end{align*} 
which concludes the proof of the theorem.
\end{proof}

\subsection{Proof for the Mixing Times}

We first recall the mixing time bound given in~\citet{atchade2011towards} using restricted spectral gaps. 

\begin{lemma}[\citet{atchade2011towards}]
\label{lem:restrictedGapConvergence}
Let \( K \) be a lazy, reversible Markov transition kernel on a state space \( \Omega \), with stationary distribution \( \Pi \). Suppose the initial distribution \( \Pi_0 \) is absolutely continuous with respect to \( \Pi \), and define 
\[
f_0(\omega) \Pi(\d \omega) = {\ \Pi_0(\d \omega)}.
\]
Assume there exist constants \( B > 1 \) and \( q > 2 \) such that \( \| f_0 \|_{\mathcal{L}^q(\Pi)} \leq B \), where \( \| \cdot \|_{\mathcal{L}^q(\Pi)} \) denotes the \( \mathcal{L}^q \)-norm with respect to \( \Pi \). Let \( \varepsilon \in (0, 1) \). Further, suppose there exists a measurable subset \( \Omega^0 \subseteq \Omega \) such that
\[
\Pi(\Omega^0) \geq 1 - \left(\frac{\varepsilon}{20B^2}\right)^{q/(q - 2)}.
\]
Then, for
\[
N \geq \frac{1}{\mathrm{SpecGap}_{\Omega^0}(K)} \log\left( \frac{2B^2}{\varepsilon^2} \right),
\]
the total variation distance between the distribution of the Markov chain \( K \) after \( N \) steps and its stationary distribution \( \Pi \) is at most \( \varepsilon \).
\end{lemma}

\begin{proof}[Proof of Lemma \ref{cor:temp}]
The first part of Equation~\eqref{eq:restricted-conv} follows directly from Lemma~\ref{lem:restrictedGapConvergence}. To prove the second part of~\eqref{eq:restricted-conv}, we first note that the TV distance between \( P^N \) and \( P \) admits the following  lower bound
\begin{align*}
\| P^N - P \|_{\text{tv}}
&= \sum_{i = 1}^{L} \int \bigl| P^N(i, \d x) - P(i, \d x) \bigr| 
\geq\int \bigl| p^N(i, x) - r_i p_i(x) \bigr| \d x, 
\quad \text{for all } i \in [L].
\end{align*}

For each $i \in [L]$, let $r_{i, N} = P^N(i, \cX)$. The TV distance between \( P_i^N \) and \( P_i \) is bounded by 
\begin{align}
\| P_i^N - P_i \|_{\text{tv}}
&= \left\| r_{i, N}^{-1} P^N(i, \cdot)  \;-\; P_i\right\|_{\text{tv}} 
= \int \left| r_{i, N}^{-1} p^N(i, x)  \;-\; p_i (x) \right| \d x 
\\ 
&\leq 
\int \left| r_{i, N}^{-1} p^N(i, x)  \;-\; r_i^{-1} p^N(i, x) \right| \d x  
+ \int \left| r_i^{-1} p^N(i, x)  \;-\; p_i (x) \right| \d x \\
&\leq 
\int \left| r_{i, N}^{-1} p^N(i, x)  \;-\; r_i^{-1} p^N(i, x) \right| \d x  + r_i^{-1} \| P^N - P \|_{\text{tv}},  
\end{align}
where the first inequality follows from the triangle inequality. 
For the first term in the last expression, using $r_{i, N} = P^N(i, \cX)$ we get 
\begin{align*}
  \int \left| r_{i, N}^{-1} p^N(i, x)  \;-\; r_i^{-1} p^N(i, x) \right| \d x 
&= r_i^{-1} \left| r_i  - r_{i, N}   \right| 
\leq \frac{r_i^{-1}}{2} \| P^N - P \|_{\text{tv}}, 
\end{align*}
where in the last step we use $r_i = P(A)$ and $r_{i, N} = P^N(A)$ with $A = \{i\} \times \cX$.  

Combining the above two displayed inequalities and using the first part of Equation~\eqref{eq:restricted-conv}, we get
\[
\| P_i^N - P_i \|_{\text{tv}} 
\;\;\leq\;\; \frac{3 }{2 r_i}\varepsilon \leq \frac{3}
     {2 \min_{ k \in [L]} r_k} \varepsilon. 
\] 
This completes the proof.    
\end{proof}
\newpage 

\section{Appendix for Section \ref{sec: example}}\label{sec:proof-example}

\subsection{Comparison of  STMH chain with approximate STMH chain} \label{sec:appx-proof-compare}
To compare the STMH chain defined in Definition~\ref{def:stmh} with the approximate STMH chain in Definition~\ref{def:approx-stmh}, it suffices to compare the stationary density $p^*_i$ with $\widetilde{p}_i$ and transition kernel $M^*_i$ with $\widetilde{M}_i$. For the former, we use Lemma 7.3 of \citet{ge2018simulated}, which shows that varying the temperature is roughly the same as changing the variance of a Gaussian distribution. 
For the latter, we derive a bound in Lemma~\ref{lm: dir}.

\begin{lemma}[Lemma 7.3 of \citet{ge2018simulated}]\label{lm: close}
Let \( 0 < \beta \leq 1 \), and suppose \( w_1, \dots, w_\m > 0 \) are weights such that \( \sum_{i = 1}^\m w_i = 1 \). Define the density functions
\[
\pi(x) \propto \left( \sum_{i = 1}^\m w_i \pi_i(x) \right)^\beta
\quad \text{and} \quad
\widetilde{\pi}(x) \propto \sum_{i = 1}^\m w_i \pi_i^\beta(x),
\]
where \( \pi_1, \dots, \pi_\m \) are component densities. Then,
\[
w_{\min} \cdot \widetilde{\pi}(x) \leq \pi(x) \leq \frac{1}{w_{\min}} \cdot \widetilde{\pi}(x),
\]
where \( w_{\min} := \min_{1 \leq i \leq \m} w_i \). 
\end{lemma}

\begin{lemma}\label{lm: dir}
 For each \( i \in [L] \), let  \( M^*_i \) be the transition kernel defined in Definition~\ref{def:stmh}, with transition density $m^*_i$. Also, let \( \widetilde{M}_i \) be the transition kernel defined in Definition~\ref{def:approx-stmh}, with transition density $\widetilde{m}_i$.
Assume that \( w_{\min} := \min_{1 \leq j \leq m} w_j > 0 \). Then, for all \( x \neq y \in \mathbb{R}^d \), the following inequality holds
\[
\widetilde{m}_i(x, y) \;\leq\; \frac{1}{w_{\min}^2} \, m^*_i(x, y).
\]
\end{lemma}

\begin{proof}
Let \( q \) denote the symmetric Gaussian proposal density  used in the Metropolis--Hastings algorithms \( M^*_i \) and \( \widetilde{M}_i \). Then, for \( x \neq y \), the transition densities are given by
\[
m^*_i(x, y) = q(x, y)\,\alpha^*_i(x, y), \quad 
\widetilde{m}_i(x, y) = q(x, y)\,\widetilde{\alpha}_i(x, y),
\]
where
\[
\alpha^*_i(x, y) = \min\!\left\{1,\; \frac{p^*_i(y)}{p^*_i(x)} \right\}, \quad 
\widetilde{\alpha}_i(x, y) = \min\!\left\{1,\; \frac{\widetilde{p}_i(y)}{\widetilde{p}_i(x)} \right\}.
\]
Using Lemma~\ref{lm: close}, we have
\[
\widetilde{p}_i(y) \leq \frac{1}{w_{\min}} p^*_i(y)
\quad \text{and} \quad
\widetilde{p}_i(x) \geq w_{\min}\, p^*_i(x),
\]
which gives
\[
\frac{\widetilde{p}_i(y)}{\widetilde{p}_i(x)} \leq \frac{1}{w_{\min}^2} \cdot \frac{p^*_i(y)}{p^*_i(x)}.
\]
Hence,
\[
\widetilde{m}_i(x, y) \leq \frac{1}{w_{\min}^2} q(x, y)\,\alpha^*_i(x, y) = \frac{1}{w_{\min}^2} m^*_i(x, y),
\]
which completes the proof.
\end{proof}

From now on, we assume that \(\cX^0 \subseteq \mathbb{R}^d \) is a measurable subset. 
Our next result, Lemma~\ref{lm: var},  shows that it suffices to obtain a lower bound on the \( [L] \times \cX^0 \)-restricted spectral gap of the approximate STMH chain in order to derive a corresponding bound for the  STMH chain. 
\begin{lemma}\label{lm: var}
Let \( \stmh \) be the STMH chain defined in Definition~\ref{def:stmh}, and let \( \widetilde{M} \) be the approximate STMH chain defined in Definition~\ref{def:approx-stmh}.
Assume that the mixture weights satisfy \( w_{\min} := \min_{1 \leq j \leq \m} w_j > 0 \). Then, the \( [L] \times \cX^0 \)-restricted spectral gaps of \( \stmh \) and \( \widetilde{M} \) satisfy the inequality
\[
\mathrm{SpecGap}_{[L] \times \cX^0}(\widetilde{M}) 
\;\leq\; \frac{1}{w_{\min}^5} \, \mathrm{SpecGap}_{[L] \times \cX^0}(\stmh).
\]
\end{lemma}

\begin{proof} 
Let \( p^* \) and \( \widetilde{p} \) denote the stationary densities of the Markov chains \( \stmh \) and \( \widetilde{M} \), respectively. Then, for all \( i \in [L] \) and \( x \in \mathbb{R}^d \), we have
\[
p^*(i, x) = r_i p^*_i(x)
\quad \text{and} \quad
\widetilde{p}(i, x) = r_i \widetilde{p}_i(x).
\]
By Lemma~\ref{lm: close}, we have 
\begin{equation}\label{eq:compare-p-s}
 w_{\min}\, \widetilde{p}_i(x) \leq  p^*_i(x) \leq w_{\min}^{-1} \,\widetilde{p}_i(x)
\end{equation}
for every $(i, x)$, and the same inequality holds for $p^*(i, x)$ and $\widetilde{p}(i, x)$ since the weights $(r_i)_{i=1}^L$ are the same for $P^*$ and $\widetilde{P}$.   Let $\mathcal{E}^*_{[L]\times\cX^0}, \widetilde{\mathcal{E}}_{[L]\times\cX^0}$ denote the ${[L]\times\cX^0}$-restricted Dirichlet forms associated with $M^*$ and $\widetilde{M}$ respectively. 
Fix a function \( g \in \mathcal{L}^2([L] \times \cX, \widetilde{p}) \)  and define \( g_i(x) := g(i, x) \) for each $(i, x)$.  
By Definition~\ref{def:restrictedSpecGap}, it suffices to show that 
\begin{equation}
   \mathrm{Var}_{P^*,{[L]\times\cX^0}}(g) \leq \frac{1}{w_{\min}^2}\mathrm{Var}_{\widetilde{P}, [L]\times\cX^0}(g), \text{ and } 
   \mathcal{E}^*_{[L]\times\cX^0}(g, g) \leq \frac{1}{w_{\min}^3} \widetilde{\mathcal{E}}_{[L]\times\cX^0}(g, g). 
\end{equation}
For the first inequality, it follows from~\eqref{eq:compare-p-s} that 
\begin{align}
\mathrm{Var}_{P^*,{[L]\times\cX^0}}(g)
&\leq \frac{1}{2 w_{\min}^2} \sum_{i=1}^{L} \sum_{j=1}^{L} \int_{\cX^0} \int_{\cX^0} \left(g(i, x) - g(j, y)\right)^2 \widetilde{p}(i, x)\, \widetilde{p}(j, y) \, \mathrm{d}x\, \mathrm{d}y\\
& =  \frac{1}{w_{\min}^2} \mathrm{Var}_{\widetilde{P},{[L]\times\cX^0}}(g).
\label{eq:firstbound}
\end{align}
By Lemma~\ref{lm: dirform}, we have 
\begin{align}\label{eq:100}
\widetilde{\mathcal{E}}_{[L]\times\cX^0}(g, g) = (1 - \lambda) \sum_{i=1}^L r_i\, \widetilde{\mathcal{E}}_{i, \cX^0}(g_i, g_i) 
  + \lambda\, \widetilde{\mathcal{E}}^I_{\cX^0}(g, g),
\end{align}
and  $\mathcal{E}^*_{{[L]\times\cX^0}}$ can be decomposed analogously. We will bound the two terms on the right-hand side of Equation~\eqref{eq:100} separately.  
For the first term, we apply Lemmas~\ref{lm: close} and ~\ref{lm: dir} to get 
\begin{align*} \widetilde{\mathcal{E}}_{i, \cX^0}(g_i, g_i)
&= \frac{1}{2}\,\int_{\cX^0} \! \int_{\cX^0}
\bigl(g_i(x) - g_i(y)\bigr)^2 \,
\widetilde{p}_i(x)\, 
\widetilde{M}_i(x, \mathrm{d}y) \, \mathrm{d}x \\
&\leq \frac{1}{2w_{\min}^3}
  \int_{\cX^0} \! \int_{\cX^0}
\bigl(g_i(x) - g_i(y)\bigr)^2 \,
p_i(x)\, 
M_i(x, \mathrm{d}y) \, \mathrm{d}x
\\&  = \frac{1}{w_{\min}^3}   \mathcal{E}_{i,\cX^0}^*(g_i, g_i). 
\end{align*} 
For the second term, we apply Lemma~\ref{lm: close} to get
\begin{align*}
\widetilde{\mathcal{E}}^I_{\cX^0}(g, g)
& = \frac{1}{4}
       \sum_{i, i' \in [L] \, : \, i' = i \pm 1} \int_{\cX^0} \bigl(g_i(x) - g_{i'}(x)\bigr)^2 \min\left\{ r_i \widetilde{p}_i(x),\, r_{i'} \widetilde{p}_{i'}(x) \right\} \, \d x\\
& \;\leq\;\frac{1}{4\,w_{\min}}
\sum_{i, i' \in [L] \, : \, i' = i \pm 1} \int_{\cX^0} \bigl(g_i(x) - g_{i'}(x)\bigr)^2 \min\left\{ r_i {p}_i(x),\, r_{i'} {p}_{i'}(x) \right\} \, \d x\\
& =  \frac{1}{w_{\min}} \, {\mathcal{E}}^{*,I}_{\cX^0}(g, g). 
\end{align*}
Combining both bounds, we obtain  that 
$\mathcal{E}^*_{[L]\times\cX^0}(g, g) \leq {w_{\min}^{-3}} \widetilde{\mathcal{E}}_{[L]\times\cX^0}(g, g)$, which concludes the proof.  
\end{proof}

\subsection{Restricted Spectral Gap of the Approximate STMH Chain} 
\label{sec:appx-proof-approx-stmh}

We begin by introducing some notation. For each \( i \in [L] \) and \( j \in [\m] \), define the \( j \)-th component of the density \( \widetilde{p}_i \) as
\begin{equation} \label{eq:tilde-pi-component}
\widetilde{p}_{(i, j)}(x) \propto \exp\left\{ -\frac{\beta_i}{2} (x - \mu_j)^\top \Sigma^{-1} (x - \mu_j) \right\},
\end{equation}
so that \( \widetilde{p}_i \) is a weighted mixture of the \( \widetilde{p}_{(i,j)} \)'s:
\[
\widetilde{p}_i(x) \propto \sum_{j=1}^{\m} w_j \widetilde{p}_{(i, j)}(x).
\]
Let \( \widetilde{M}_{(i,j)} \) denote the Metropolis--Hastings transition kernel targeting \( \widetilde{p}_{(i, j)} \), with a symmetric Gaussian proposal  density \( q(x, y) = \mathcal{N}(y; x, \eta I) \), where \( \eta > 0 \) is the step size.
  We set
\[
\cX^0 := \left\{ x \in \mathbb{R}^d \;:\; \|x\| \leq R \right\},
\]
where \( R > 0 \) is a fixed radius. To obtain a lower bound on the \([L] \times \cX^0 \)-restricted spectral gap of the approximate STMH chain, we invoke Theorem~\ref{thm: ddt}. 
Assumption~\ref{ass:pi} holds by our construction of $\widetilde{P}$. The following lemmas verify the other two assumptions required for this theorem. 

\subsubsection{Validation of Condition (\ref{eq:decomp-c1}) in Assumption~\ref{ass:Mij}
}\label{sec:appx-proof-c1}

\begin{lemma}\label{lm: l1}
For each \( i \in [L] \), let \( \widetilde{p}_i \) be the density defined in Equation~\eqref{eq:tilde-pi}, and let  \( g_i \in \mathcal{L}^2(\cX, \widetilde{p}_i) \). Then the following inequality holds
\begin{equation*}  
    \sum_{j=1}^\m w_j\,\widetilde{\mathcal{E}}_{(i,j), \cX^0}(g_i, g_i),
    \;\leq\;\widetilde{\mathcal{E}}_{i, \cX^0}(g_i, g_i) \qquad \forall i \in [L],
\end{equation*}
where \( \widetilde{\mathcal{E}}_{(i,j), \cX^0} \) denotes the $\cX^0$-restricted Dirichlet form of the kernel \( \widetilde{M}_{(i,j)} \), and \( \widetilde{\mathcal{E}}_{i, \cX^0} \) denotes the $\cX^0$-restricted Dirichlet form of the kernel \( \widetilde{M}_i \), as defined in Definition~\ref{def:approx-stmh}. 

In particular, for the approximate STMH chain \( \widetilde{M} \) defined in Definition~\ref{def:approx-stmh}, condition~\eqref{eq:decomp-c1} holds  with constant \( C_1 = 1 \). 
\end{lemma}
\begin{proof}
For any nonnegative real sequences \(\{a_j\}\) and \(\{b_j\}\), we have the inequality

$\min\left\{
\sum_{j} a_j,\;\sum_{j} b_j
\right\}
\;\ge\;
\sum_{j} \min\left\{a_j,\,b_j\right\}.$

Applying this  to $\widetilde{p}_i= \sum_{j = 1}^\m w_j \widetilde{p}_{(i, j)}$, we obtain
\begin{align}
\min\left\{ \widetilde{p}_i(x), \widetilde{p}_i(z) \right\}
\;\geq\;
\sum_{j = 1}^\m w_j \min\left\{ \widetilde{p}_{(i, j)}(x), \widetilde{p}_{(i, j)}(z) \right\},
 \label{eq:random}
\end{align}
for all  $x, z \in \mathbb{R}^d$  and  $i \in [L]$. Let \( q \) denote the symmetric Gaussian proposal density associated with the Metropolis--Hastings algorithms \( \widetilde{M}_i \) and \( \widetilde{M}_{(i,j)} \). Then, applying Equation~\eqref{eq:random}, we obtain
\[
\begin{aligned}
\widetilde{\mathcal{E}}_{i, \cX^0}(g_i,g_i) 
&=
\frac{1}{2} 
\int_{\cX^0} \int_{\cX^0}
\left(g_i(x)-g_i(z)\right)^2
q(x,  z)\min\left\{\widetilde{p}_i(x), \widetilde{p}_i(z)\right\} \, \d x \, \d z\\
&\geq
\frac{1}{2} 
\int_{\cX^0} \int_{\cX^0}
\left(g_i(x)-g_i(z)\right)^2 q(x, z)
\sum_{j=1}^\m w_j \min\left\{\widetilde{p}_{(i, j)}(x), \widetilde{p}_{(i, j)}(z)\right\} \, \d x \,  \d z\\
&=
\sum_{j=1}^\m w_j \widetilde{\mathcal{E}}_{(i,j), \cX^0}(g_i,g_i).
\end{aligned}
\]
This completes the proof of the Lemma.
\end{proof}

\subsubsection{Validation of Condition (\ref{eq:poincare-Mij}) in Assumption~\ref{ass:Mij}
}\label{sec:appx-proof-c2}

We lower bound the \( \cX^0 \)-restricted spectral gap of each Metropolis--Hastings chain \( \widetilde{M}_{(i,j)} \) using the path method of~\citet{yuen2000applications} in the following lemma.

\begin{lemma}\label{lm: l2} Let \( 0 < \eta \leq R^2  \).
For each \( i \in [L] \) and \( j \in [\m] \), the Markov chain \( \widetilde{M}_{(i,j)} \) admits the following lower bound on its \( \cX^0 \)-restricted spectral gap 
\[
\mathrm{SpecGap}_{\cX^0}\bigl(\widetilde{M}_{(i,j)}\bigr)
\;\geq\; \frac{\gamma_{\min}^{d/2} \eta^{3/2}}{ 13 R^{d+3}} \, .
\]

In particular, for the approximate STMH chain \( \widetilde{M} \) defined in Definition~\ref{def:approx-stmh}, condition~\eqref{eq:poincare-Mij} holds with constant $$C_2 = \frac{ 13 R^{d+3}} {\gamma_{\min}^{d/2}}. $$
\end{lemma}

\begin{proof}
We use the linear path method described in Section~2 of~\cite{yuen2000applications}.  
This approach also extends to the restricted spectral gap setting; see, for example,~\citet{atchade2021approximate} and~\citet{chang2024dimension} where the canonical path method has been adapted to the restricted spectral gap in discrete spaces. For each pair \( (x, y) \in \cX^0 \times \cX^0 \), we construct a linear path connecting \( x \) to \( y \), with all intermediate points lying in \( \cX^0 \). Fix a step size \( \delta > 0 \), and define the number of steps along the path by
\[
b_{xy} \coloneqq \left\lceil \frac{\|x - y\|}{\delta} \right\rceil.
\]
The path is then given by
\[
\gamma_{xy} = \left( \gamma_{xy}^{(0)}, \dots, \gamma_{xy}^{(b_{xy})} \right),
\]
where
\[
\gamma_{xy}^{(i)} \coloneqq \frac{(b_{xy} - i)x + i y}{b_{xy}}, \; \text{ for } 0 \le i \leq b_{xy}.
\] 
Let \(\Gamma \coloneqq \{ \gamma_{xy} : (x, y) \in \cX^0 \times \cX^0\}\) denote the collection of all such paths, and let \(E\) denote the set of all edges that appear in at least one path \(\gamma_{xy} \in \Gamma\). 
The capacity of an edge \( (u, v) \in E \) is given by
\[
T_{(i, j)}(u, v) \coloneqq \widetilde{p}_{(i,j)}(u) \, \widetilde{m}_{(i,j)}(u, v)
= \widetilde{p}_{(i,j)}(u) \, q(u, v) \min\!\left\{1,\;\frac{\widetilde{p}_{(i,j)}(v)}{\widetilde{p}_{(i,j)}(u)}\right\},
\]
where \( \widetilde{m}_{(i,j)} \) denotes the transition density corresponding to the kernel \( \widetilde{M}_{(i,j)} \), and \( q(u, v) \) is the Gaussian proposal density associated with kernel $\widetilde{M}_{(i, j)}$. As shown in Section~2 of~\cite{yuen2000applications}, the set of paths \( \Gamma \) satisfies the regularity conditions and, for any \( (u, v) \in \gamma_{xy} \), the associated Jacobian satisfies
$
J_{x,y}(u, v) = b_{xy}^d$ {(see~\cite[page 5]{yuen2000applications} for details)}. Then, by Theorem 2.1 and Corollary 2.4 in~\cite{yuen2000applications}, we have 
\begin{equation}\label{eq:spec}
\mathrm{SpecGap}_{\cX^0}\bigl(\widetilde{M}_{(i,j)}\bigr) \geq \frac{1}{A} \end{equation} where 
\[A = \operatorname*{ess\,sup}_{(u,v) \in E} 
\left\{
\frac{1}{T_{(i, j)}(u,v)} 
\sum_{\gamma_{xy} \ni (u,v)}  |\gamma_{xy}|\,
\widetilde{p}_{(i,j)}(x) \, \widetilde{p}_{(i,j)}(y) \, b_{xy}^d 
\right\},
\] 
and $|\gamma_{xy}|$ denotes the length of the path $\gamma_{xy}$. Since  $\tilde{p}_{(i, j)}$ is log-concave, for any $(u, v) \in \gamma_{xy}$, we have  
\[
\min\left\{ \widetilde{p}_{(i,j)}(x), \, \widetilde{p}_{(i,j)}(y) \right\}
\leq 
\min\left\{ \widetilde{p}_{(i,j)}(u), \,\widetilde{p}_{(i,j)}(v) \right\}. 
\]
Hence, $T_{(i, j)}(u, v) \geq q(u, v) \min \{ \widetilde{p}_{(i,j)}(x), \, \widetilde{p}_{(i,j)}(y)  \}$, and we can upper bound \(A\) as
\begin{align}
A \leq 
b^{d+1} \cdot \operatorname*{ess\,sup}_{(u,v) \in E} 
\left\{ q(u, v)^{-1}
\sum_{\gamma_{x,y} \ni (u,v)} 
\widetilde{p}_{(i,j)}(z_{x,y})
\right\}, \label{eq: K}
\end{align}
where
$
b \coloneqq \max_{(x, y) \in \cX^0 \times \cX^0} b_{x,y} 
$ and \(z_{x,y}\) is defined as
\[
z_{x,y} :=
\begin{cases}
x, & \text{if } \max\{ \widetilde{p}_{(i,j)}(x), \widetilde{p}_{(i,j)}(y) \} = \widetilde{p}_{(i,j)}(x), \\
y, & \text{otherwise}.
\end{cases}
\]
Note that
\[
\widetilde{p}_{(i,j)}(z_{x,y}) 
\;\leq\; 
\frac{\beta_i^{d/2}}{(2\pi \gamma_{\min}
)^{d/2}},
\] and the proposal density $q(u, v)$ is given by
\[
q(u,v) 
\;=\; 
\frac{1}{(2\pi \eta)^{d/2}} \,
\exp\!\left(-\frac{\|v - u\|^2}{2\eta}\right).
\]
Substituting this into Equation \eqref{eq: K}, we obtain
\[
A \;\leq\; 
\frac{ \beta_i^{d/2} \eta^{d/2} b^{d + 3}}{\gamma_{\min}^{d/2}} \cdot 
\operatorname*{ess\,sup}_{(u,v) \in E} 
\left\{ 
\exp\left(\frac{\|v - u\|^2}{2\eta}\right)
\right\},  
\]
where we have also used that an edge $(u, v)$ belongs to at most $b^2$ paths in $\Gamma$.  
Choose step size $\delta=   \sqrt{5 \eta}$, which yields
\begin{equation}
    b \leq \Big\lceil \frac{2R}{  \sqrt{5 \eta}} \Big\rceil  \leq \frac{R}{\sqrt{\eta}}. 
\end{equation}
Since $\beta_i \leq 1$,  we obtain that 
\[
A \;\leq\;  \frac{\beta_i^{d/2}R^{d+3}}{\gamma_{\min}^{d/2} \eta^{3/2}} e^{2.5} \leq 
  \frac{ 13 R^{d+3}}{\gamma_{\min}^{d/2} \eta^{3/2}}.
\]
From Equation~\eqref{eq:spec}, we get
\[
\mathrm{SpecGap}_{\cX^0}\bigl(\widetilde{M}_{(i,j)}\bigr) \;\geq\; \frac{1}{A}
\;\geq\; \frac{\gamma_{\min}^{d/2} \eta^{3/2}}{ 13 R^{d+3}},
\]
which concludes the proof.
\end{proof}

\subsubsection{Auxiliary Lemmas}

To verify Assumption~\ref{ass:barM} and compute the constant $C_3$ in condition~\eqref{eq:poincare-barM}, 
we will need several lemmas. The proof of Lemma~\ref{lm:can} is omitted. 

\begin{lemma}[Canonical Paths Bound]\label{lm:can}
Let \( \mathcal{S} \) be a finite state space, and let \( K \) be the transition kernel of a reversible Markov chain on \( \mathcal{S} \) with stationary distribution \( \Pi \) and Dirichlet form \( \mathcal{E} \). For each pair of distinct states \( x, y \in \mathcal{S} \), let \( \gamma_{xy} \) denote a path from \( x \) to \( y \) consisting of valid transitions under \( K \), i.e.,
\[
x = x_0 \to x_1 \to x_2 \to \dots \to x_{n-1} \to x_n = y.
\]
Let \( \Gamma = \{ \gamma_{xy} : x, y \in \mathcal{S},\ x \neq y \} \) be the collection of such paths for all distinct pairs \( (x, y) \). The edge congestion associated with \( \Gamma \) is defined as
\[
\rho_e(\Gamma) = \max_{\substack{u,v \in \mathcal{S} \\ K(u,v) > 0}} \frac{1}{\Pi(u) K(u,v)} \sum_{\substack{(u,v) \in \gamma_{xy} \\ \gamma_{xy} \in \Gamma}} \Pi(x) \Pi(y) |\gamma_{xy}|,
\]
where $|\gamma_{xy}|$ denotes the length of the path $\gamma_{xy}$. Then, for any function \( g : \mathcal{S} \to \mathbb{R} \), the following Poincar\'{e} inequality holds
\[
\operatorname{Var}_{\Pi}(g) \leq \rho_e(\Gamma) \, \mathcal{E}(g, g).
\]
\end{lemma}

\begin{lemma}\label{lm:min-TV}
Let $\Pi, \tilde{\Pi}$ be two probability distributions (absolutely continuous with respect to each other) with density function $\pi, \tilde{\pi}$ respectively. Then, 
\begin{equation}
    \int  \min \left\{ \pi(x),  \tilde{\pi}(x) \right\}  \d x 
    = 1 - \frac{1}{2} \| \Pi - \tilde{\Pi} \|_{\rm{tv}}
    \geq 1 -  \sqrt{ \frac{1}{2} \mathrm{KL}( \Pi \mid \tilde{\Pi} )  }.  
\end{equation}
\end{lemma}

\begin{lemma}  \label{lm:KL-normal}
Let \(|\Sigma|\) denote the determinant of a matrix \(\Sigma\). 
The Kullback-Leiber divergence between two $d$-dimensional Gaussian distributions with equal means is given by 
\begin{equation}
\mathrm{KL} \left(  N(\mu, \Sigma_1) \mid  N(\mu, \Sigma_2) \right)  = \frac{1}{2} \left\{ \log \frac{ |\Sigma_2| }{|\Sigma_1|} - d + \mathrm{tr}(\Sigma_2^{-1} \Sigma_1) \right\}. 
\end{equation} 
\end{lemma}   
  
\begin{lemma}\label{lm:H}
Let
\(
D \coloneqq \max\left\{ \max_{k \in [\m]} \|\mu_k\|,\, \sqrt{\gamma_{\min}} \right\}.
\)
For each \( i \in [L] \), \( j \in [\m] \), define
\[
\widetilde{p}_{i}(j \mid x) := \frac{w_j \, \widetilde{p}_{(i,j)}(x)}{\sum\limits_{k=1}^\m w_k \, \widetilde{p}_{(i,k)}(x)} \qquad \text{for all } x \in \mathbb{R}^d,
\]
where \( \widetilde{p}_{(i,j)}(x) \) denotes the density defined in Equation~\eqref{eq:tilde-pi-component}. Then, for all \( i \in [L] \), \( j \in [\m] \), and \( x \in \mathbb{R}^d \), the following inequalities hold
\begin{align}
\label{eq:Z21}
\widetilde{p}_{(i,j)}(x) 
&\geq \left( \frac{\beta_i}{2\pi \gamma_{\max}} \right)^{d/2} 
\exp\left( - \frac{\beta_i}{2 \gamma_{\min}} ( \|x\| + D )^2 \right), \\
\label{eq:Z}
\widetilde{p}_i(j \mid x) 
&\geq w_j \exp\left( - \frac{\beta_i}{\gamma_{\min}} ( \|x\| + D )^2 \right).
\end{align}
\end{lemma}
\begin{proof}
 To prove Equation~\eqref{eq:Z21}, we write
\[
\widetilde{p}_{(i,j)}(x) \geq \left( \frac{\beta_i}{2\pi\gamma_{\max}} \right)^{d/2} \exp\left( - \frac{\beta_i}{2} (x - \mu_j)^\top \Sigma^{-1} (x - \mu_j) \right),
\]
where the inequality follows from the bound \( |\Sigma| \leq \gamma_{\max}^d \). Next, we use the inequality
\[
\|(x - \mu_j)^\top \Sigma^{-1} (x - \mu_j) \|
\leq \frac{1}{\gamma_{\min}} \|x - \mu_j\|^2 
\leq \frac{1}{\gamma_{\min}} (\|x\| + \|\mu_j\|)^2 
\leq \frac{1}{\gamma_{\min}} (\|x\| + D)^2 
\]
to obtain
\[
\widetilde{p}_{(i,j)}(x) 
\geq \left( \frac{\beta_i}{2\pi\gamma_{\max}} \right)^{d/2} 
\exp\left( - \frac{\beta_i}{2 \gamma_{\min}} (\|x\| + D)^2 \right).
\]
This establishes Equation~\eqref{eq:Z21}. 
To prove Equation~\eqref{eq:Z}, we define the function \( \widetilde{J} : \mathbb{R}^d \to [\m] \) by
\[
\widetilde{J}(x) := \arg\max_{k \in [\m]} \, \widetilde{p}_{(i,k)}(x).
\]It follows that for every \( k \in [\m] \),
\[
\widetilde{p}_{(i,k)}(x) \leq \widetilde{p}_{(i,\widetilde{J}(x))}(x),
\]
and therefore,
\[
\sum_{k=1}^\m w_k \, \widetilde{p}_{(i,k)}(x) \leq \sum_{k=1}^\m w_k \, \widetilde{p}_{(i,\widetilde{J}(x))}(x) = \widetilde{p}_{(i,\widetilde{J}(x))}(x).
\]
Substituting this upper bound into the  definition of \( \widetilde{p}_i(j \mid x)\), we get
\begin{equation}\label{eq:3}
\widetilde{p}_{i}(j \mid x) \geq  \frac{w_j\widetilde{p}_{(i,j)}(x)}{\widetilde{p}_{(i,\widetilde{J}(x))}(x)} .
\end{equation}
To simplify the ratio of Gaussian densities on the right-hand side, we expand the expression explicitly as
\begin{align*} 
 \frac{\widetilde{p}_{(i,j)}(x)}{\widetilde{p}_{(i,\widetilde{J}(x))}(x)}  = \exp\left\{
- \beta_i(\mu_{\widetilde{J}(x)} - \mu_{j})^\top \Sigma^{-1} x
- \frac{\beta_i}{2} \left( \mu_{j}^\top \Sigma^{-1} \mu_{j} - \mu_{\widetilde{J}(x)}^\top \Sigma^{-1} \mu_{\widetilde{J}(x)} \right)
\right\}.
\end{align*}
By definition of $D$, we have
$
\| \mu_{\widetilde{J}(x)} - \mu_{j} \| \leq 2 D.
$ Using the Cauchy--Schwarz inequality, we obtain
\[
\| (\mu_{\widetilde{J}(x)} - \mu_{j})^\top \Sigma^{-1} x \| 
\leq \| \mu_{\widetilde{J}(x)} - \mu_{j} \| \cdot \| \Sigma^{-1} x \| 
\leq \frac{2D}{\gamma_{\min}} \| x \|,
\]
and similarly,
\[
\left\| \mu_{j}^\top \Sigma^{-1} \mu_{j} - \mu_{\widetilde{J}(x)}^\top \Sigma^{-1} \mu_{\widetilde{J}(x)}\right\|
\leq \frac{1}{\gamma_{\min}} \left( \| \mu_{j} \|^2 + \| \mu_{\widetilde{J}(x)} \|^2 \right)
\leq \frac{2 D^2}{\gamma_{\min}}.
\]
Putting these together, we get
\begin{equation}\label{eq:2}
  \frac{\widetilde{p}_{(i,j)}(x)}{\widetilde{p}_{(i,\widetilde{J}(x))}(x)} 
\geq \exp\left(-\frac{\beta_i}{\gamma_{\min}} ( 2D \|x\| + D^2 ) \right) \geq \exp\left( - \frac{\beta_i}{\gamma_{\min}} (  \|x\| + D )^2 \right).  
\end{equation}

Equations~\eqref{eq:3} and \eqref{eq:2} together prove  Equation~\eqref{eq:Z}. This completes the proof.
\end{proof}

\subsubsection{Validation of Condition (\ref{eq:poincare-barM}) in Assumption~\ref{ass:barM}
}\label{sec:appx-proof-c3}

Let \( \widehat{M} \) denote the projected chain associated with the approximate STMH chain \( \widetilde{M} \). We next establish a lower bound on the spectral gap of \( \widehat{M} \) using the canonical paths  method.  
Recall that we define $\cX^0 = \left\{ x \in \mathbb{R}^d  :  \|x\| \leq R \right\}$.

\begin{lemma}\label{lm:ddg}
Let \( \widetilde{M} \) denote the approximate STMH chain defined in Definition~\ref{def:approx-stmh}.  
Define the following parameters
\[
D \coloneqq \max\left\{ \max_{k \in [\m]} \|\mu_k\|,\, \sqrt{\gamma_{\min}} \right\}, \quad r \coloneqq \frac{\min_{i \in [L]} r_i}{\max_{i \in [L]} r_i}.
\]
Suppose the following conditions hold
\begin{itemize}
  \item[(i)] Let $R \geq \sqrt{d}D$ be such that $P_{(i, j)}(\cX^0) \geq {3}/{4}$  for all  $i \in [L]$ and $j \in [\m]$,
\item[(ii)] \( \beta_1 = \Theta\left( {\gamma_{\min}}/{D^2} \right) \) and \( \beta_1 \leq 1 \),
  \item[(iii)] ${\beta_{i+1}}/{\beta_i} \leq 1 + {1}/{\sqrt{d}}$ \quad for all $i \in [L - 1]$.
\end{itemize}
 Let \( \widehat{M} \) be the projected chain defined in Definition~\ref{def:project-chain} associated with    $\widetilde{M}$. Under these conditions, \( \widehat{M} \)  satisfies the spectral gap bound
\[
\mathrm{SpecGap}(\widehat{M}) 
\;\geq\; 
\frac{3\min\{(1 - \lambda), \lambda\}\,r^2}{64{L}^2 \kappa^{d/2} \exp(cd)},
\]
where $c >0$ is a fixed constant. In particular, for the approximate STMH chain \( \widetilde{M} \), condition~\eqref{eq:poincare-barM} holds with constant
\[
C_3 = \frac{64{L}^2 \kappa^{d/2} \exp(cd)}{3\min\{(1 - \lambda), \lambda\}\,r^2}.
\]
\end{lemma}

\begin{proof}
We construct the canonical paths as follows. 
Fix two arbitrary states \( x = (i, j) \), \( y = (i', j') \in [L]\times [\m]\) with $i \leq i'$. 
\setlength{\leftmargini}{1cm}
\begin{itemize} 
    \item[(a)] If \( j = j' \), let $\gamma_{x  y}$ be  
    \( (i, j) \rightarrow  (i+1, j) \rightarrow \ldots \rightarrow (i', j) \).  
    \item[(b)] If \( j \neq j' \), let $\gamma_{x  y}$ be \( (i, j) \rightarrow (i-1, j) \rightarrow \ldots \rightarrow (1, j) \rightarrow (1, j') \rightarrow (2, j') \rightarrow \ldots \rightarrow (i', j') \).
\end{itemize}
Define $\gamma_{y x}$ as the reverse of $\gamma_{x y}$.
Let \( \Gamma \) denote the collection of such paths over all distinct pairs \( (x, y)   \).  Let \( i \in [L] \), \( j \in [n] \), and \( i' = i \pm 1 \in [L] \). From the definition of \( \widehat{M} \), we have
\begin{align*}
    \widehat{M}((i, j), (i', j)) 
    &= \frac{\lambda}{2} \int_{\cX^0} \frac{ \widetilde{p}_{(i,j)}(x) }{ \widetilde{P}_{(i, j)} (\cX^0) } \cdot \widetilde{a}\big( (i, j, x), (i', j, x) \big) \, \mathrm{d}x,
\end{align*}
where the acceptance probability is given by
\[
\widetilde{a}\big( (i, j, x), (i', j, x) \big) = \min \left\{ \frac{r_{i'} \, \widetilde{p}_{(i', j)}(x)}{r_i \, \widetilde{p}_{(i, j)}(x)}, \; 1 \right\}.
\]
Hence, the probability of transitioning from state \( (i, j) \) to \( (i - 1, j) \) under the projected chain \( \widehat{M} \) is given by \begin{align*}
    \widehat{M}((i, j), (i - 1, j)) 
    &= \frac{\lambda}{2 \widetilde{P}_{(i, j)}(\cX^0)} \int_{\cX^0} \min\left\{ \frac{r_{i - 1}}{r_i} \, \widetilde{p}_{(i - 1,j)}(x), \, \widetilde{p}_{(i,j)}(x) \right\} \, \mathrm{d}x,
\end{align*}
for all  $i \in \{2, \dots, L\}$ and $ j \in [n]$. Since \( {r_{i-1}}/{r_i} \geq r \) by definition and $\widetilde{P}_{(i, j)}(\cX^0)\leq 1 $, we have
\begin{equation}\label{eq:lowerbound}
\widehat{M}((i, j), (i - 1, j)) 
\geq \frac{\lambda r}{2} \int_{\cX^0} \min\left\{ \widetilde{p}_{(i - 1,j)}(x), \, \widetilde{p}_{(i,j)}(x) \right\} \, \mathrm{d}x.
\end{equation}
By Lemma~\ref{lm:min-TV} and Lemma~\ref{lm:KL-normal}, 
\begin{align*}
    \int_{\mathbb{R}^d} \min\left\{ \widetilde{p}_{(i - 1,j)}(x), \, \widetilde{p}_{(i,j)}(x) \right\} \, \mathrm{d}x  &\geq 1 - \sqrt{ \frac{1}{2} \mathrm{KL} (\widetilde{P}_{(i,j)} \mid  \widetilde{P}_{(i-1,j)})  }  \\
    &\geq 1 - \frac{\sqrt{d}}{2}     \sqrt{ f\left(  \frac{\beta_i}{\beta_{i-1}}\right) }, 
    \text{ where } f(x) = x - 1 - \log x.    
\end{align*}
For $x \geq 1$, we have $f(x) \leq (x-1)^2/2$. Hence,   if $\beta_i / \beta_{i-1} - 1 = 1 / \sqrt{d}$, then 
\begin{align*}
     \int_{\mathbb{R}^d} \min\left\{ \widetilde{p}_{(i - 1,j)}(x), \, \widetilde{p}_{(i,j)}(x) \right\} \, \mathrm{d}x  \geq \frac{1}{2}.
\end{align*} 
Condition (i) implies 
 \begin{align*}
     \int_{\cX^0} \min\left\{ \widetilde{p}_{(i - 1,j)}(x), \, \widetilde{p}_{(i,j)}(x) \right\} \, \mathrm{d}x  \geq \frac{1}{4}.
\end{align*} 
 Substituting the above bound  in Equation~\eqref{eq:lowerbound}, we obtain
\[
\widehat{M}((i, j), (i - 1, j)) 
\geq \frac{\lambda r}{8}.
\]

 Similarly, we can derive the bound
\[
\widehat{M}((i, j), (i + 1, j)) \geq \frac{\lambda r}{8}, \qquad \text{for all } i \in [L - 1],\; j \in [n]. \]
Next, we derive a lower bound on the transition probability from \( (1, j) \) to \( (1, j') \) in the projected chain \( \widehat{M} \), which is given by
\begin{align*}
\widehat{M}((1, j), (1, j')) = 
(1 - \lambda) \int_{\cX^0}\frac{ \widetilde{p}_{(1,j)}(x) }{ \widetilde{P}_{(1,j)} (\cX^0) } \cdot 
\widetilde{p}_{1}(j' \mid x)
\, \mathrm{d}x, \qquad j, j'\in[\m].
\end{align*}
By applying Lemma~\ref{lm:H} and noting that $\widetilde{P}_{(1, j)}(\cX^0) \leq 1$, we obtain
\begin{align*}
\widehat{M}((1, j), (1, j')) &\geq 
(1 - \lambda)w_{j'} \left( \frac{\beta_1}{2\pi \gamma_{\max}} \right)^{d/2} \int_{\cX^0}  \exp\left( - \frac{2\beta_1}{\gamma_{\min}} ( \|x\| + D)^2 \right)
\, \mathrm{d}x.
\end{align*}
By condition~(ii), there exist fixed constants \( c_1, c_2 > 0 \) such that
\[
c_1\, \frac{\gamma_{\min}}{D^2} < \beta_1 < c_2\,\frac{ \gamma_{\min}}{D^2}.
\]
Let 
$
\cX^D \coloneqq \left\{ x = (x_1, \dots, x_d) \in \mathbb{R}^d : |x_i| \leq D \text{ for all } i \in [d] \right\} \subseteq \cX^0.$ Then for any \( x \in \cX^D \), we have 
\(
\|x\| + D \leq 2\sqrt{d}D,
\)
which implies
\[
\exp\left( - \frac{2\beta_1}{\gamma_{\min}} (\|x\| + D)^2 \right) 
\geq \exp\left( - \frac{8\beta_1 d D^2}{\gamma_{\min}} \right) \ge \exp(-8c_2d).
\] Therefore, we obtain the following lower bound\begin{align*}
\widehat{M}((1, j), (1, j')) 
&\geq (1 - \lambda) w_{j'} (c_1)^{d/2} \left( \frac{\gamma_{\min}}{2\pi \gamma_{\max} D^2} \right)^{d/2} 
\exp(-8c_2d) \cdot \mathrm{Vol}(\cX^D),
\end{align*}
where \( \mathrm{Vol}(\cX^D) \) denotes the volume of the cube \( \cX^D \). Substituting \( \mathrm{Vol}(\cX^D) = (2D)^d \), we get
\begin{align}
\widehat{M}((1, j), (1, j')) 
&\geq  (1 - \lambda) w_{j'} (c_1)^{d/2} \left( \frac{\gamma_{\min}}{2\pi \gamma_{\max} D^2} \right)^{d/2} 
\exp(-8c_2d) \cdot (2D)^d \\
&\geq (1 - \lambda) w_{j'} \cdot \kappa^{-d/2} \cdot \exp(-c d), \label{eq:lowbound}
\end{align}
where  \( c > 0 \) is a fixed constant. 

Let \( \gamma_{xy} \) be a path between any two vertices \( x, y \in [L]\times [\m]\). Then, \( |\gamma_{xy}| \leq 2L \).
We now derive an upper bound on the edge congestion \( \rho_e(\Gamma) \) defined in Lemma~\ref{lm:can}. Let $z, w \in [L] \times [n]$. 
\begin{itemize}
 \item[(a)] 
 Let \( z = (i,j) \) and \( w = (i - 1,j) \). Then the edge \( (z,w) \) is used only by paths between vertices  \( x\) and \( y \) such that one lies in the set
  $
  S := \{i, \ldots, L\} \times \{j\}, \text{and the other in } S^c.
  $
  Therefore, its contribution to the edge congestion is bounded by
  \[
  \frac{
    \sum_{(x,y) \in \Gamma: ((i,j), (i - 1,j)) \in \gamma_{xy}} 
    |\gamma_{xy}|\,\widehat{P}(x)\, \widehat{P}(y)
  }{
    \widehat{P}((i,j))\, \widehat{M}((i,j), (i - 1,j))
  }
  \leq
  \frac{(2L)\, \widehat{P}(S)\, \widehat{P}(S^c)}{
    \widehat{P}((i,j))\, \widehat{M}((i,j), (i - 1,j))
  },
  \]
  where \[
\widehat{P}((\ell, k)) = r_\ell w_k \frac{\widetilde{P}_{(\ell, k)}(\cX^0)}{\widetilde{P}([L] \times \cX^0)}, \qquad \ell \in [L], \; k \in [n]
\]
denotes the stationary distribution of the projected chain \( \widehat{M} \).
 We have the following bounds
\[
\frac{\widehat{P}(S)}{\widehat{P}((i,j))} 
= \frac{\widehat{P}(\{i, \ldots, L\} \times \{j\})}{\widehat{P}((i,j))} 
\leq \frac{4L}{3r}, \qquad 
\widehat{P}(S^c) \leq 1, \qquad 
\widehat{M}((i,j), (i - 1,j)) \geq \frac{\lambda r}{8},
\]
where the first bound follows from condition (i).
 Combining these, we conclude
  \[
  \frac{
    \sum_{(x,y) \in \Gamma: ((i,j), (i - 1,j)) \in \gamma_{xy}} 
    |\gamma_{xy}|\, \widehat{P}(x)\, \widehat{P}(y)
  }{
    \widehat{P}((i,j))\, \widehat{M}((i,j), (i - 1,j))
  }
  \leq \frac{64L^2}{3\lambda r^2}.
  \]
Similarly, we obtain the same bound for \( z = (i, j) \) and \( w = (i + 1, j) \).
\item[(b)] Let \( z = (1,j) \) and \( w = (1,j') \). Then the edge \( (z,w) \) is used only by paths between vertices $x$ and $y$ such that one of them lies in the set \( [L] \times \{j\} \) and the other in \( [L] \times \{j'\} \). Therefore, its contribution to edge congestion is bounded by
\[
\frac{
  \sum_{(x,y) \in \Gamma: ((1,j), (1,j')) \in \gamma_{xy}} 
  |\gamma_{x,y}|\, \widehat{P}(x)\, \widehat{P}(y)
}{
  \widehat{P}((1,j))\,\widehat{M}((1,j), (1,j'))
}
\leq
\frac{
  (2L)\, \widehat{P}([L] \times \{j\})\,
         \widehat{P}([L] \times \{j'\})
}{
  \widehat{P}((1,j))\, \widehat{M}((1,j), (1,j'))
}.
\]
We now bound each term on the right-hand side. By condition~(i), we have
\begin{align*}
\widehat{P}([L] \times \{j\}) 
&\leq \frac{4L}{3r} \, \widehat{P}((1,j)), 
\qquad
\widehat{P}([L] \times \{j'\}) \leq \frac{4 w_{j'}}{3},
\end{align*}
and from Equation~\eqref{eq:lowbound}, we have
\begin{align*}
\widehat{M}((1,j), (1,j')) 
&\geq (1 - \lambda) w_{j'} \cdot \kappa^{-d/2} \cdot \exp(-cd).
\end{align*}
Combining these, we obtain
\[
\frac{
  \sum_{(x,y): ((1,j),(1,j')) \in \gamma_{xy}} 
  |\gamma_{xy}|\,\widehat{P}(x)\, \widehat{P}(y)
}{
 \widehat{P}((1,j))\, \widehat{M}((1,j), (1,j'))
}
\leq  \frac{32{L}^2 \kappa^{d/2} \exp(cd)}{9(1 - \lambda)r} .
\]
\end{itemize}
Let \( \lambda \in (0, 1) \) be a fixed constant.
 Thus, the edge congestion associated with $\Gamma$ is bounded by $$\rho_e(\Gamma) \leq \frac{64{L}^2 \kappa^{d/2} \exp(cd)}{3\min\{(1 - \lambda), \lambda\}\,r^2}.$$ From Lemma \ref{lm:can}, projected chain $\widehat{M}$  satisfies the Poincar\'{e} inequality 
\[
\operatorname{Var}_{\widehat{P}}(\widehat{g}) \leq \frac{64{L}^2 \kappa^{d/2} \exp(cd)}{3\min\{(1 - \lambda), \lambda\}\,r^2}\,\mathcal{\widehat{E}}(\widehat{g}, \widehat{g}),  
\quad \forall \, \widehat{g} \colon [L] \times [\m] \rightarrow \mathbb{R}, 
\]
where \( \widehat{\mathcal{E}} \) denotes the Dirichlet form associated with the projected chain \( \widehat{M} \). This completes the proof of the lemma. 
\end{proof}

\subsubsection{Restricted Spectral Gap Bound}

We now invoke Theorem~\ref{thm: ddt} to bound the \([L] \times \cX^0\)-restricted spectral gap of the approximate STMH chain $\widetilde{M}$, as formalized in the next lemma.

\begin{lemma}\label{lm:aproxpoin}
Let \( \widetilde{M} \) denote the approximate STMH chain, as defined in Definition~\ref{def:approx-stmh}, with $\lambda$ being a fixed constant. 
Under the same conditions as in Lemma~\ref{lm:ddg}, the \( [L] \times \cX^0 \)-restricted spectral gap of \( \widetilde{M} \) admits the following lower bound

\[
\mathrm{SpecGap}_{[L] \times \cX^0}(\widetilde{M}) 
\;\geq\; 
\Omega\left(\frac{\gamma_{\min}^{d/2}r^2 \eta^{3/2}}{R^{d + 3}{{L}^2 \kappa^{d/2} \exp(cd)}}\right),
\]
where $c >0$ is a fixed constant. 
\end{lemma}

\begin{proof}
For the approximate STMH chain \( \widetilde{M} \), Assumption~\ref{ass:Mij} is satisfied with constants \( C_1 = 1 \) and \[ C_2 = \frac{ 13 R^{d+3}}{\gamma_{\min}^{d/2} \eta^{3/2}}.\]Additionally, Assumption~\ref{ass:barM} holds with \[ C_3 =  \frac{64{L}^2 \kappa^{d/2} \exp(cd)}{3\min\{(1 - \lambda), \lambda\}\,r^2}, \] 
where $c>0$ is a fixed constant. Since $\lambda$ is treated as fixed, combining these with Theorem~\ref{thm: ddt} completes the proof of the lemma.
\end{proof}

\subsection{Restricted Spectral Gap of the STMH chain}
\label{sec:appx-proof-stmh}

We now establish a lower bound on the \( [L] \times \cX^0 \)-restricted spectral gap of the STMH chain \( \stmh, \) as formalized in next lemma.

\begin{lemma} \label{thm: poin}
Let  \( \stmh \) be the STMH chain, as defined in Definition~\ref{def:stmh}, with $\lambda$ being a fixed constant. 
Under the same conditions as in Lemma~\ref{lm:ddg}, the \( [L] \times \cX^0 \)-restricted spectral gap of \( \stmh \) admits the lower bound
\[
\mathrm{SpecGap}_{[L] \times \cX^0}(\stmh) 
\;\geq\; 
\Omega\left(\frac{w_{\min}^5\gamma_{\min}^{d/2}r^2 \eta^{3/2}}{R^{d + 3}{{L}^2 \kappa^{d/2} \exp(cd)}}\right),
\]
where $c >0$ is a fixed constant. 
\end{lemma}
\begin{proof}
This follows directly from Lemma~\ref{lm: var} and Lemma~\ref{lm:aproxpoin}.
\end{proof}

\subsection{Estimation of Partition Functions}
\label{sec:appx-proof-partition}

\subsubsection{Assumptions on the Parameters} \label{sec:ass_para}
We now describe how to choose the algorithm parameters---number of temperature levels $L$, inverse temperature sequence $(\beta_i)_{i=1}^L$, temperature-swap rate $\lambda$, proposal step size $\eta$, initial covariance matrix $\Sigma_0$, number of iterations $N$---so that the STMH algorithm achieves the asserted time complexity.  
Recall that $\kappa = \gamma_{\max} / \gamma_{\min}$. 

\begin{align}
  L  &= \Theta\left[ 
        \kappa 
        \left\{ 
            D^2 + \log w_{\min}^{-1}
            + d \left( 1 + \log \kappa \right) 
        \right\} 
        \log\left( \frac{D^2}{\gamma_{\min}} \right) + 1 
     \right], \label{eq:L} \\
\beta_1 &= \Theta\left(\frac{\gamma_{\min}}{D^2}\right),  
\hspace{0.1cm} 
\frac{\beta_{i+1}}{\beta_i}  \leq \min\left\{
1 + \frac{1}{\sqrt{d}}, \;
\frac{ \gamma_{\min}}{D^2 + 2 \gamma_{\max} d \nu } 
\right\}
\hspace{0.1cm}\text{for } i \in [L - 1], \label{eq:b}  \\
& \quad\quad\quad \text{ where } \nu = 1 + \log \kappa + \frac{2}{d} \log \left(\frac{2}{w_{\min}}\right), \\ 
  R &= D + 
     \sqrt{d \kappa D^2} + 
     \sqrt{2 \kappa D^2 \log \left( \frac{20\, e^6 L^2\kappa^d}{w_{\min}^2\varepsilon} \right)}, \label{eq:R} \\
     N &\geq 
    \frac{C' L^4 R^d\kappa^{d/2}\exp(c'd)}{\gamma_{\min}^{d/2} w_{\min}^5} 
    \log\left( 
      \frac{L^2\kappa^d}{\varepsilon^2 w_{\min}^2} 
    \right), \hspace{0.1cm}\text{for some fixed constants } c', C' > 0
, \label{eq:N} \\
    \Sigma_0 &= \sigma_0^2 I, 
\hspace{0.1cm}
  \text{where } \sigma_0^2 =\Theta\left(\frac{\gamma_{\min}}{\beta_1}\right), 
    \label{eq:S} \\
\lambda &\text{ is any fixed constant}, \label{eq:lambda} \\
\eta  &\geq R^2. \label{eq:eta}
\end{align}

\subsubsection{Auxiliary Lemmas} 
\begin{lemma}
\label{lem:ri-bounds}
Let  $L > 0$ be an integer. Assume the partition–function estimates $\widehat Z_1,\dots,\widehat Z_L$ satisfy
\begin{equation}
 \frac{\widehat Z_i / Z_i}{\widehat Z_1 / Z_1} 
\in\Bigl[(1-\tfrac1L)^{i-1},\,(1+\tfrac1L)^{i-1}\Bigr],
\qquad \text{for all } i \in [L].
 \label{eq:est-bias}
\end{equation}
Define
\[
r_i \;:=\;
\frac{Z_i/\widehat Z_i}{\displaystyle\sum_{k=1}^{L} Z_k/\widehat Z_k},
\qquad \text{for all } i \in [L].
\]
Then, 
\begin{equation}
\frac{e^{-2}}{L}\;\le\;r_i\;\le\;\frac{e^{2}}{L},
\qquad \text{for all } i \in [L]. \label{eq:ri-unif}
\end{equation}
Moreover, \[\displaystyle r := \frac{\min_{i \in [L]} r_i}{\max_{i \in [L]} r_i} \ge e^{-4}.\]

\end{lemma}

\begin{proof}
For each \( i \in [L] \), define
\(
b_i := {Z_i}/{\widehat{Z}_i}
\),
and denote their sum by \( S \),
\[
S := \sum_{k=1}^{L} b_k.
\]
Then \( r_i = b_i / S \). From Equation~\eqref{eq:est-bias} we have, for every $i \in [L]$,
\begin{equation}
(1+\tfrac1L)^{-(i-1)}
\;\le\;
\frac{b_i}{b_1}
\;\le\;
(1-\tfrac1L)^{-(i-1)},
\label{eq:bi-ratio}
\end{equation}
which gives
\begin{equation}
L\,b_1\,(1+\tfrac1L)^{-(L-1)}
\;\le\;
S
\;\le\;
L\,b_1\,(1-\tfrac1L)^{-(L-1)}.
\label{eq:SumBounds}
\end{equation}
Combining Equations~\eqref{eq:bi-ratio} with~\eqref{eq:SumBounds}, we get
\[
r_i
=\frac{b_i}{S}
\;\ge\;
\frac{b_1(L+1)^{-(L-1)}}
     {L\,b_1\,(L-1)^{-(L-1)}}
= \frac1L\,
   \Bigl(\frac{L-1}{L+1}\Bigr)^{L-1},
\]
and 
\[
r_i
\;\le\;
\frac{b_1(L-1)^{-(L-1)}}
     {L\,b_1\,(L+1)^{-(L-1)}}
= \frac1L\,
   \Bigl(\frac{L+1}{L-1}\Bigr)^{L-1}.
\]
Define 
\[
C_L \;:=\,\Bigl(\frac{L+1}{L-1}\Bigr)^{L-1}.
\]
Taking logarithms and using $\log(1+x)\le x$ for all $x>-1$, we obtain
\[
\log C_L = (L-1)\,\log\,\!\Bigl(1+\tfrac{2}{\,L-1}\Bigr)
          \;\le\; (L-1)\Bigl(\tfrac{2}{\,L-1}\Bigr)
          \;\leq\;2,
\]
so $C_L \le e^{2}$.  Likewise $C_L^{-1}\ge e^{-2}$.  Therefore
\[
\frac{e^{-2}}{L}\;\le\;r_i\;\le\;\frac{e^{2}}{L},
\qquad \text{for all } i \in [L],
\]
establishing~\eqref{eq:ri-unif}. This further implies that \( r \geq e^{-4} \), completing the proof of the lemma.
\end{proof}

\begin{lemma}\label{lem:L2-f0-tilde}
Let \( L > 0 \) be an integer, and suppose \( \beta_1 \) and \( \sigma_0 \) satisfy Equation~\eqref{eq:b} and Equation~\eqref{eq:S}, respectively.  Then the initial density is given by
\[
p^0(1, \cdot) = \mathcal{N}\left(0, \sigma_0^2I\right),
\qquad
p^0(i, \cdot) = 0 \quad \text{for all } i \in [L] \setminus \{1\},
\]
with the corresponding distribution denoted by \( P^0 \). The stationary density of the STMH chain \( \stmh \), as defined in Definition~\ref{def:stmh}, is
\[
p(i, x) = r_i\, p_i^*(x), \quad  i \in [L],\, x \in \mathbb{R}^d,
\]
where  \( (r_i)_{i = 1}^L \) are defined in Lemma~\ref{lem:ri-bounds}, and the component densities \( (p_i^*)_{i \in [L]} \) are given in Equation~\eqref{eq:pi}. Let \( P \) denote the corresponding stationary distribution. Define 
\(
f_0 := {dP^0}/{dP}.
\)
Then,
\[
\|f_0\|_{\mathcal{L}^2(P)}^2 \leq \frac{c_1 \,\exp(c_2d)\,L\,\kappa^{d/2}}{w_{\min}},
\]
for some fixed constants $c_1, c_2>0.$
\end{lemma}

\begin{proof}
The \( \mathcal{L}^2 \)-norm of \( f_0 \) is given by
\begin{equation}\label{eq:f0}
\|f_0\|_{\mathcal{L}^2(P)}^2
= \sum_{i=1}^{L} \int p(i, x) \, |f_0(i, x)|^2 \, dx
= \frac{1}{r_1} \int \frac{(p^0(1, x))^2}{p_1^*(x)} \, dx,
\end{equation}
where the second equality follows from the fact that \( P^0 \) is supported only on  \( i = 1 \). By Lemma~\ref{lm: close}, for all \( x \in \mathbb{R}^d \), we have
\begin{equation}\label{eq:Z_1}
p_1^*(x)
\;\ge\;
w_{\min} \, \widetilde{p}_1(x) = w_{\min} \sum_{j = 1}^{\m} w_j \widetilde{p}_{(1, j)}(x).
\end{equation}
where $\widetilde{p}_1$ is defined in Equation~\eqref{eq:tilde-pi} and $\widetilde{p}_{(1, j)}$ is defined in Equation~\eqref{eq:tilde-pi-component}.
Substituting \eqref{eq:Z_1} into \eqref{eq:f0}, we obtain
\[
\|f_0\|_{\mathcal{L}^2(P)}^2
\leq \frac{1}{r_1 w_{\min}} \int \frac{(p^0(1, x))^2}{\sum_{j = 1}^{\m} w_j \widetilde{p}_{(1, j)}(x)} \, dx.
\]
Using the convexity of the \( \chi^2 \)-divergence, we further bound this as
\begin{equation}\label{eq:boundf0}
\|f_0\|_{\mathcal{L}^2(P)}^2
\leq \frac{1}{r_1 w_{\min}} \sum_{j = 1}^{\m} w_j \int \frac{(p^0(1, x))^2}{\widetilde{p}_{(1, j)}(x)} \, dx.
\end{equation}

For each \( j \in [\m] \), the density \( \widetilde{p}_{(1, j)}(x) \) can be lower bounded as
\begin{align}\label{eq:x}
\widetilde{p}_{(1, j)}(x)
\geq
\left( \frac{\beta_1}{2\pi\gamma_{\max}} \right)^{d/2}
\exp\left( -\frac{\beta_1}{2\gamma_{\min}} \|x - \mu_j\|^2 \right).
\end{align}
Since \( p^0(1, \cdot) \sim \mathcal{N}(0, \sigma_0^2 I) \), we have
\[
(p^0(1, x))^2
=
\left( \frac{1}{2\pi\sigma_0^2 } \right)^d 
\exp\left( -\frac{1}{\sigma_0^2 } \|x\|^2 \right).
\]
From Equation~\eqref{eq:S}, we have fixed constants $0 <s_1, s_2 <2$ such that 
\[
s_1 \frac{\gamma_{\min}}{\beta_1} \leq\sigma_0^2 \leq s_2 \frac{\gamma_{\min}}{\beta_1}.
\]
This gives 
\begin{align}\label{eq:y}
(p^0(1, x))^2
\leq
\left( \frac{\beta_1}{2\pi s_1 \gamma_{\min}} \right)^d 
\exp\left( -\frac{\beta_1}{s_2 \gamma_{\min}} \|x\|^2 \right).
\end{align}
Substituting Equations~\eqref{eq:x} and~\eqref{eq:y} into Equation~\eqref{eq:boundf0}, we obtain 
\begin{align*}
&\|f_0\|_{\mathcal{L}^2(P)}^2\\
&\leq \frac{\kappa^{d/2}}{s_1^d\,r_1 w_{\min}} 
\sum_{j = 1}^{\m} w_j 
\left( \frac{\beta_1}{2\pi\gamma_{\min}} \right)^{d/2} 
\exp\left( \frac{\beta_1}{(2 - s_2)\gamma_{\min}} \|\mu_j\|^2 \right)
\int 
\exp\left( -\frac{\beta_1(2 - s_2)}{2\gamma_{\min}s_2} \left\lVert x + \frac{s_2\,\mu_j}{2 - s_2}\right\rVert^2 \right) \, \d x \\
&= \frac{\kappa^{d/2} \,s_2^{d/2}}{s_1^d\, (2 - s_2)^{d/2} \,r_1 w_{\min}} 
\sum_{j = 1}^{\m} w_j 
\exp\left( \frac{\beta_1}{(2 - s_2)\gamma_{\min}} \|\mu_j\|^2 \right).
\end{align*}
From Equation~\eqref{eq:b}, we have a fixed constant $s_3> 0$ such that \( \beta_1 \leq s_3{\gamma_{\min}}/{D^2} \). Substituting this, we get
\[
\|f_0\|_{\mathcal{L}^2(P)}^2 \leq
\frac{\kappa^{d/2} \,s_2^{d/2}}{s_1^d\, (2 - s_2)^{d/2} \,r_1 w_{\min}} \exp\left(\frac{s_3}{2 - s_2}\right).
\]
 By Lemma~\ref{lem:ri-bounds}, we have $r_1 \geq 1/(e^2L)$. Substituting this into the above bound on \( \|f_0\|_{\mathcal{L}^2(P)}^2 \) proves the lemma.
\end{proof}
\begin{lemma}
Let $X$ follow the $d$-dimensional Gaussian distribution  with mean $\mu$ and covariance matrix $\Sigma$. Denote the largest eigenvalue of $\Sigma$ by $\|\Sigma\|$. Then, 
\begin{align*}
    \mathbb{P} \left(  \| X  \| \leq \| \mu \| + 
    \sqrt{d \|\Sigma\| } + \sqrt{ 2 \|\Sigma\| \log (1 / \varepsilon) } \right) \geq \varepsilon. 
\end{align*}
\end{lemma}
\begin{proof}
Using the standard concentration inequality for Lipschitz functions of Gaussian random vectors, we get 
\begin{equation}
     \mathbb{P} \left(  \| X - \mu \| \geq \mathbb{E}( \| X - \mu \|) + t \right) \leq e^{- t^2/ 2 \|\Sigma\|}.  
\end{equation}
Since $\mathbb{E}( \| X - \mu \|) \leq  \sqrt{d \|\Sigma\| }$, we get 
\begin{equation}
     \mathbb{P} \left(  \| X - \mu \| \geq \sqrt{d \|\Sigma\| } + t \right) \leq e^{- t^2/ 2  \|\Sigma\|}.  
\end{equation}
Letting $t = \sqrt{ 2 \gamma_{\mathrm{max}} \log (1 / \varepsilon) }$ and applying triangle inequality, we get the asserted bound. 
\end{proof}
\begin{lemma}\label{lm:initial}
Suppose \( 1 \leq \ell \leq L \), and let Algorithm~\ref{alg:simulated_tempering} be run with the potential function \( f(x) \) defined in Equation~\eqref{eq:potential}, inverse temperatures $\beta_1\ < \cdots < \beta_\ell$, and using the parameters specified in Equations~\eqref{eq:b}, \eqref{eq:R}, \eqref{eq:N}, \eqref{eq:S}, \eqref{eq:lambda} and \eqref{eq:eta}. Assume that the partition function estimates \(\widehat{Z}_1, \dots, \widehat{Z}_\ell\) satisfy Equation~\eqref{eq:est-bias}. Let \( P^N \) denote the distribution obtained after running Algorithm~\ref{alg:simulated_tempering}  for \( N \) steps, and let \( P \) denote its stationary distribution. Then the total variation distance between \( P \) and \( P^N \) satisfies
\begin{equation}\label{eq:tvbound}
\| P - P^N \|_{\mathrm{tv}} \;\le\; \varepsilon.
\end{equation}
\end{lemma}
\begin{proof}
Under the assumptions of the lemma, and by Lemmas~\ref{thm: poin} and~\ref{lem:ri-bounds}, 
the \([\ell] \times \cX^0\)-restricted spectral gap of $\stmh$ satisfies
\[
\mathrm{SpecGap}_{[\ell] \times \cX^0}(\stmh) 
\;\geq\; 
\Omega\left( \frac{w_{\min}^5 \, \gamma_{\min}^{d/2}}{R^d \, \ell^4 \, \kappa^{d/2} \, \exp(cd)} \right),
\]
where $c > 0$ is a fixed constant. Moreover, from Lemma~\ref{lem:L2-f0-tilde}, $\|f_0\|_{\mathcal{L}^2(P)}^2$ is bounded above by $B$, where
\[
B = \frac{c_1\, \exp(c_2d) \ell \kappa^{d/2}}{w_{\min}},
\]
and $c_1, c_2>0$ are  fixed constants. Applying Lemma~\ref{lem:restrictedGapConvergence} with the above parameters yields the desired total variation bound.
\end{proof} 
\begin{lemma}\label{lm: partition}
Assume the same conditions and notations as in Lemma~\ref{lm:initial}. Let \( P^*_\ell \) denote the marginal stationary distribution at temperature level \( \ell \), with density 
\(
p^*_\ell(x) \propto \exp(-\beta_\ell f(x)),
\)
and let \( P^N_\ell \) denote the marginal distribution at level \( \ell \) after running Algorithm~\ref{alg:simulated_tempering} for \( N \) steps, with density
\(
p^N_\ell(x) \propto P^N(\ell, x).
\)
Then the total variation distance between \( P^*_\ell \) and \( P^N_\ell \) is bounded by
\[
\| P^*_\ell - P^N_\ell \|_{\mathrm{tv}} \;\le\; \frac{3e^2\ell}{2}\, \varepsilon.
\]
\end{lemma}
\begin{proof}
By Lemma~\ref{lem:ri-bounds}, we have \( \min_{i \in [\ell]} r_i \geq 1 / (e^2 \ell) \). The proof now follows directly from Lemmas~\ref{cor:temp} and~\ref{lm:initial}.
\end{proof}
 {In the following lemma, we analyze how many times Algorithm~\ref{alg:simulated_tempering} must be re-run, with a fixed number of steps \( N \), in order to obtain a sample from the desired temperature level.}

\begin{lemma}
\label{lem:sample-complexity-temperature}
Suppose the partition function estimates \(\widehat{Z}_1, \dots, \widehat{Z}_\ell\) satisfy Equation~\eqref{eq:est-bias}. Let \( I_N \in [\ell] \) denote the temperature index of the state returned after running Algorithm~\ref{alg:simulated_tempering} for \(N\) steps. Suppose the algorithm is run independently \(T\) times, each for \(N\) steps. Then, for any fixed temperature level \(k \in [\ell]\), if
\[
T \;\ge\; e^{2} \ell \log\!\left(\tfrac{1}{\alpha}\right), \qquad \alpha \in (0,1),
\]
the probability that at least one of the \(T\) runs returns a sample from level \(k\) satisfies
\[
{\P}\left( \exists\, t \in [T] \text{ such that } I^{(t)}_N = k \right) \;\ge\; 1 - \alpha,
\]
where \(I^{(t)}_N\) is the temperature level returned in the \(t\)-th run.
\end{lemma}
\begin{proof}
Let  \(k \in [\ell]\). From Lemma~\ref{lem:ri-bounds}, we have
\[
{\P}(I_N \neq k) \;=\; 1 - {\P}(I_N = k)
\;\le\; 1 - \frac{1}{e^{2}\,\ell}.
\]
Hence,

\[
{\P}\left( \not\exists\, t \in [T] \text{ such that } I^{(t)}_N = k \right)
\leq \bigl(1 - \tfrac{1}{e^{2}\,\ell}\bigr)^{T}
\;\le\;
\exp\!\bigl(-\tfrac{T}{e^{2}\,\ell}\bigr).
\]
Setting this upper bound no larger than \( \delta \), and solving for \( T \), completes the proof of the lemma.
\end{proof}

Assuming the partition function estimates satisfy Equation~\eqref{eq:est-bias}, we have shown that the algorithm reaches total variation distance at most \( \varepsilon \) within the time complexity specified in Equation~\eqref{eq:N}. We now show that partition function estimates satisfy Equation~\eqref{eq:est-bias}. By combining these two components, we establish the overall time complexity for the complete algorithm.

\begin{lemma}\label{lm:estimation}
Let \( \delta \in (0, 1) \) and $1\leq \ell\leq L.$ Suppose the parameters satisfy Equations~\eqref{eq:L}, \eqref{eq:b}, \eqref{eq:R},  \eqref{eq:S}, \eqref{eq:lambda}, and \eqref{eq:eta}, and assume that the partition function estimates \( \widehat{Z}_1, \dots, \widehat{Z}_\ell \) satisfy Equation~\eqref{eq:est-bias}. Let \( s = L^2 \log(1/\delta) \). Collect \( s \) samples from Algorithm~\ref{alg:simulated_tempering}, denoted by \( (x_j)_{j = 1}^s \). Define the next partition function estimate \( \widehat{Z}_{\ell + 1} \) as
\[
\widehat{Z}_{\ell+1} := \overline{r}\,\widehat{Z}_\ell, \quad \text{where} \quad \overline{r} := \frac{1}{s} \sum_{j=1}^s
\exp\bigl(-(\beta_{\ell+1} - \beta_\ell) f(x_j) \bigr).
\]
Then, with probability at least \( 1 - \delta \), the estimate \( \widehat{Z}_{\ell+1} \) also satisfies Equation~\eqref{eq:est-bias}. In particular,
\begin{equation}
\left| \frac{\widehat{Z}_{\ell + 1} / Z_{\ell + 1}}{\widehat{Z}_1 / Z_1} \right|
\in
\left[
  \left(1 - \frac{1}{L}\right)^{\ell},
  \left(1 + \frac{1}{L}\right)^{\ell}
\right].
\end{equation}
\end{lemma}
The proof of Lemma \ref{lm:estimation} requires the following results.
\begin{lemma}[Lemma 9.1 of \citet{ge2018simulated}] \label{lem:9.1}
Suppose that \(P_1\) and \(P_2\) are probability measures on \(\Omega\) with density functions (with respect to a reference measure)
\[
p_1(x) = \frac{g_1(x)}{Z_1}, \quad \text{and} \quad  p_2(x) = \frac{g_2(x)}{Z_2}.
\]
Suppose \(\widetilde{P}_1\) is a measure such that \(\|\widetilde{P}_1 - P_1\|_{\text{tv}} < {c}/{2C^2}\), and \({g_2(x)}/{g_1(x)} \in [0, C]\) for all \(x \in \Omega\).
Given \(n\) samples \(x_1, \dots, x_n\) from \(\widetilde{P}_1\), define the random variable
\[
\overline{r} = \frac{1}{n} \sum_{i=1}^n \frac{g_2(x_i)}{g_1(x_i)}. 
\]
Let
\[
r = \E_{x \sim P_1} \frac{g_2(x)}{g_1(x)} = \frac{Z_2}{Z_1}. 
\]
and suppose \(r \geq {1}/{C}\). Then with probability at least \(1 - e^{-{n c^2}/{(2C^4)}}\),
\[
\left| \frac{\overline{r}}{r} - 1 \right| \leq c. 
\]
\end{lemma}

\begin{lemma}[Lemma G.16 of \citet{ge2018simulated}] \label{lem:G.15}
Suppose that 
\(
f(x) = -\log \left[ \sum_{i=1}^n w_i \, e^{-f_i(x)} \right],
\)
where \( f_i(x) = f_0(x - \mu_i) \), and \( f_0 \colon \mathbb{R}^d \to \mathbb{R} \) is a \(\kappa\)-strongly convex and \(K\)-smooth function. For any \( a > 0 \), let \( P_a \) denote the probability measure with density
\(
p_a(x) \propto e^{-a f(x)}.
\)
Let \( Z_a \) be the corresponding normalization constant, given by
\(
Z_a = \int_{\mathbb{R}^d} e^{-a f(x)} \, \mathrm{d}x.
\) Suppose that \( \|\mu_i\| \leq D \) for all \( i \in [n] \), and let \( \alpha, \beta > 0 \). Let $$ A = D + \frac{1}{\sqrt{\alpha \kappa}} 
\left( \sqrt{d} + \sqrt{d \log \left(\frac{K}{\kappa}\right) + 2 \log \left(\frac{2}{w_{\min}}\right)} \right).$$
If $ \alpha < \beta,$ then 
$$
 \min_{x\in \mathbb{R}^d} \frac{p_\alpha(x)}{p_\beta(x)} \geq \frac{Z_\beta}{Z_\alpha} \qquad \text{and} \qquad \frac{Z_\beta}{Z_\alpha}  \in \left[ \frac{1}{2} e^{-\frac{1}{2} (\beta - \alpha) K A^2}, 1 \right]. $$
\end{lemma}
\begin{proof}[Proof of Lemma \ref{lm:estimation}] By Equation~\eqref{eq:b} and Lemma~\ref{lem:G.15}, we have
\[
\frac{\exp(-\beta_{\ell + 1} f(x))}{\exp(-\beta_{\ell} f(x))} = \exp\left(-(\beta_{\ell+1} - \beta_\ell) f(x)\right) \in [0,1/(2\,e)]
\]
for all \( \ell \in [L - 1] \). Moreover, by substituting \( \varepsilon = {4}/({3\ell L}) \) into Lemma~\ref{lm: partition}, we obtain

\[
\| P^*_\ell - P^{\widetilde{N}}_\ell \|_{\mathrm{tv}} \;\le\; \frac{2e^2}{L},
\]
when $$\widetilde{N} \geq 
    \frac{C'L^4 R^d\kappa^{d/2}\exp(c'd)}{\gamma_{\min}^{d/2} w_{\min}^5} 
    \log\left( 
      \frac{L^4\kappa^d}{w_{\min}^2} 
    \right), $$ 
    where $C', c' >0$ are fixed constants.  Next, by applying Lemma~\ref{lem:9.1} with constants \( C ={1}/{2e} \) and \( c ={1}/{L} \), we obtain the following bound 
\[
\left| \frac{\widehat{Z}_{\ell + 1} / Z_{\ell + 1}}{\widehat{Z}_\ell / Z_\ell} \right|
\in
\left[
  1 - \frac{1}{L},\;
  1 + \frac{1}{L}
\right].
\]

The lemma then follows by induction on \( \ell \).
\end{proof}
\subsubsection{Proof of Theorem \ref{thm: pa}}\label{sec:main_thm}
\begin{proof}[{Proof of Theorem \ref{thm: pa}}]
Let \( L \) denote the number of temperature levels defined in Equation~\eqref{eq:L}. By applying Lemma~\ref{lm:estimation} inductively with \( \delta = \varepsilon / (4L) \), we obtain that, with probability at least \( 1 - \varepsilon / 4 \), the following bound holds
\[
\frac{\widehat{Z}_{\ell}}{Z_\ell} \in \left[ \left( 1 - \frac{1}{L} \right)^{\ell - 1}, \left( 1 + \frac{1}{L} \right)^{\ell - 1} \right] \cdot \frac{\widehat{Z}_{1}}{Z_1}
\qquad \text{for all } \ell \in [L].
\]
To ensure this guarantee, it suffices to generate
\(
s = L^2 \log\left( {4L}/{\varepsilon} \right)
\)
samples from each temperature level \( i \in [L] \), resulting in a total of
\(
sL = L^3 \log\left( {4L}/{\varepsilon} \right)
\)
samples from Algorithm~\ref{alg:simulated_tempering}.  Applying Lemma~\ref{lem:sample-complexity-temperature} with \( \alpha = {\varepsilon}/{(4L^4 \log(4L / \varepsilon))} \), we obtain that, with probability at least \( 1 - \varepsilon / 4 \), we obtain \( s \) samples from each temperature level $i \in [L]$ by running Algorithm~\ref{alg:simulated_tempering} for \( N \) steps (as defined in Equation~\eqref{eq:N}) and repeating this process independently $T$ times, where
\[
T = sL \cdot e^2 L \log\left( \frac{1}{\alpha} \right)
= e^2 L^4 \log\left( \frac{4L}{\varepsilon} \right) \log\left( \frac{4L^4}{\varepsilon} \log\left(\frac{4L}{\varepsilon}\right) \right).
\]
Hence, the total time complexity for getting partition function estimates is
\[
T_{\text{partition}} = T \cdot N
=   \frac{C'\,L^8 R^d\kappa^{d/2}\exp(c'd)}{\gamma_{\min}^{d/2} w_{\min}^5} 
    \log^3\left( 
      \frac{L\kappa}{\varepsilon w_{\min}} 
    \right),
\]
where $c', C' >0$ are fixed constants. 
By applying Lemma~\ref{lm: partition} and Lemma~\ref{lem:sample-complexity-temperature}, we conclude that, with probability at least \( 1 - \varepsilon/4 \), Algorithm~\ref{alg:simulated_tempering} produces a sample from a distribution that is within total variation distance \( \varepsilon/4 \) of the target distribution \( P^* \) in time \( T_{\text{sampling}} \), where
\begin{align}
  T_{\text{sampling}} &= 
  e^2 \,L \, \log\left(\frac{4}{\varepsilon}\right)\frac{C''L^4 R^d\kappa^{d/2}\exp(c''d)}{\gamma_{\min}^{d/2} w_{\min}^5} 
    \log\left( \frac{L^2\kappa^d}{\varepsilon^2 w_{\min}^2} \right)\\
    &= \frac{C'L^5R^d\kappa^{d/2}\exp(c'd)}{\gamma_{\min}^{d/2} w_{\min}^5} 
    \log^2\left( \frac{L\kappa}{\varepsilon w_{\min}} \right),
  \end{align}
where $c', C', c'', C'' >0$ are fixed constants. The overall time complexity \( T \) consists of two components: the time to get partition function estimates, and the time to generate sample from the target distribution 
\[
T = T_{\text{partition}} + T_{\text{sampling}}.
\]
This completes the proof of the theorem.
 
\end{proof}

\end{document}